\documentclass[10pt]{amsart}
\usepackage[australian]{babel}
\usepackage[
    letterpaper,
    top=3.25cm,
    bottom=2cm,
    left=3cm,
    right=3cm,
    marginparwidth=1.75cm
]{geometry}

\usepackage[T1]{fontenc}
\usepackage{lmodern}
\usepackage{graphicx}
\usepackage{amsthm,amsmath,amssymb}
\usepackage{mathtools}
\usepackage{microtype}
\usepackage{setspace}
\usepackage{xcolor}
\usepackage{enumerate}
\usepackage{booktabs,longtable,array}
\usepackage{listings}
\usepackage{xurl}

\usepackage{tikz}

\usepackage[
    colorlinks=true,
    citecolor=black,
    linkcolor=black,
    urlcolor=blue
]{hyperref}

\theoremstyle{plain}
\newtheorem{theorem}{Theorem}[section]
\newtheorem{proposition}[theorem]{Proposition}
\newtheorem{lemma}[theorem]{Lemma}
\newtheorem{corollary}[theorem]{Corollary}

\theoremstyle{remark}
\newtheorem{remark}[theorem]{Remark}

\newtheoremstyle{main}
  {}{}{\itshape}{}{\bfseries}{.}{5pt plus 1pt minus 1pt}
  {\thmname{#1}\thmnote{ {\normalfont(#3)}}}
\theoremstyle{main}
\newtheorem{maintheorem}{Main Theorem}

\theoremstyle{plain}

\newcommand{\maintheoremautorefname}[1]{Main Theorem}
\AtBeginDocument{%
}

\newcommand{\nc}{\newcommand}

\nc{\RR}{\mathbb{R}}
\nc{\CC}{\mathbb{C}}
\nc{\ZZ}{\mathbb{Z}}
\nc{\CP}{\mathbb{CP}}

\DeclareMathOperator{\Ric}{Ric}
\DeclareMathOperator{\Hess}{Hess}
\DeclareMathOperator{\tr}{tr}
\DeclareMathOperator{\diag}{diag}
\DeclareMathOperator{\Ad}{Ad}
\DeclareMathOperator{\Id}{Id}

\nc{\SU}{\mathsf{SU}}
\nc{\PSU}{\mathsf{PSU}}
\nc{\dd}{d}
\nc{\norm}[1]{\left\lVert #1\right\rVert}
\nc{\abs}[1]{\left\lvert #1\right\rvert}
\nc{\Fc}{{\mathcal F_{\mathrm{loc}}}}

\title[A compact non-K\"ahler shrinking Ricci soliton]
{A compact non-K\"ahler gradient shrinking Ricci soliton}
\author[Timothy Buttsworth and Henry Wood]
{Timothy Buttsworth$^1$ and Henry Wood$^2$}

\makeatletter
\g@addto@macro\@setauthors{%
  \begingroup
  \centering\normalfont\footnotesize
  \noindent\textsuperscript{1}School of Mathematics and Statistics,
  The University of New South Wales, Kensington, Sydney, NSW 2052, Australia.\\
  E-mail address: \href{mailto:t.buttsworth@unsw.edu.au}{\nolinkurl{t.buttsworth@unsw.edu.au}}\par
  \smallskip
  \noindent\textsuperscript{2}School of Mathematics and Physics,
  The University of Queensland, St Lucia, Brisbane, QLD 4072, Australia.\\
  E-mail address: \href{mailto:henry.wood@uq.edu.au}{\nolinkurl{henry.wood@uq.edu.au}}\par
  \endgroup
  \medskip
}
\makeatother

\subjclass[2020]{Primary 53E20; Secondary 53C25, 65G20}
\keywords{Ricci soliton, cohomogeneity one, flag manifold, computer-assisted proof}

\begin{document}

\begin{abstract}
We construct a gradient shrinking Ricci soliton on
the simply-connected eight-manifold
$\SU(3)\times_{\mathsf T^2}\CP^1$. This Riemannian manifold is
not Einstein, is nowhere locally K\"ahler, is locally irreducible, 
and is invariant under the effective cohomogeneity one action of $\PSU(3)\times S^1$. 
This new gradient Ricci soliton is found by first constructing an approximate solution, certifying the residual, and applying a fixed-point theorem.  
\end{abstract}

\maketitle
\setcounter{tocdepth}{1}
\tableofcontents

\section{Introduction}\label{sec:introduction}

A \emph{Ricci soliton} on a smooth manifold is a pair \((g,X)\) consisting of a Riemannian metric \(g\) and a vector field \(X\) satisfying
\begin{equation}\label{eq:ricci-soliton}
  \Ric_g+\frac12\mathcal L_Xg=\lambda g,
\end{equation}
for some constant \(\lambda \in \RR\). 
Under Ricci flow, Ricci solitons evolve only by scaling and pullback by diffeomorphisms and thus provide models for singularities of the Ricci flow.
A Ricci soliton is called \emph{shrinking}, \emph{steady}, or \emph{expanding} if the constant \(\lambda\) in \eqref{eq:ricci-soliton} is positive, zero, or negative.
If \(X\) is a gradient field of some smooth scalar function $u$, then the Ricci soliton equation \eqref{eq:ricci-soliton} becomes
\begin{equation}\label{eq:gradient-ricci-soliton-eq}
  \Ric_g+\Hess_g u=\lambda g.
\end{equation}
In this case, we call the metric \(g\) a \emph{gradient Ricci soliton} with \emph{soliton potential} \(u\).
On a compact Ricci soliton, \(X\) is gradient up to addition of a Killing field. If the Ricci soliton is non-Einstein, it is shrinking \cite[\S\S2--3]{Perelman}.

Many examples of compact non-Einstein Ricci solitons have been
constructed using K\"ahler metrics \cite{Koiso,Cao,DancerWang}.
Tian and Cao asked whether there exists a compact Ricci soliton that is nowhere locally K\"ahler, non-Einstein, and locally irreducible \cite[Question~9.1]{Derdzinski}.
Here \emph{nowhere locally K\"ahler} means that no nonempty open subset
carries a compatible parallel complex structure, and \emph{locally
irreducible} means that no nonempty open subset is isometric to a product
of two positive-dimensional Riemannian manifolds. 
Let $\mathsf T^2\subset\SU(3)$ act on $\CP^1$ by \(\diag(z_1,z_2,z_3)\cdot[a:b]=[a:(z_3/z_1)b]\), and define
\begin{equation}\label{eq:manifold-intro}
  \qquad
  M=\SU(3)\times_{\mathsf T^2}\CP^1.
\end{equation}

\begin{maintheorem}\label{thm:main}
The compact simply-connected manifold $M$ admits a smooth gradient
shrinking Ricci soliton whose metric is nowhere locally
K\"ahler, non-Einstein, and locally irreducible.
\end{maintheorem}

The manifold $M$ already carries a K\"ahler shrinking Ricci soliton
from the construction of Dancer and Wang \cite{DancerWang}.
Their work extends the projective-line bundle constructions of
Koiso \cite{Koiso} and Cao \cite{Cao} to bundles over products of Fano
K\"ahler--Einstein manifolds and generalised flag manifolds.
We construct a non-K\"ahler shrinking soliton on \(M\)
by solving \eqref{eq:gradient-ricci-soliton-eq} for metrics and
potentials invariant under the cohomogeneity one action of
$\SU(3)\times S^1$ on $M$. This reduces the construction to a shooting problem for a system
of ordinary differential equations, with smoothness conditions
at the two singular orbits. 
B\"ohm \cite{Bohm} used this approach to construct infinitely many
inhomogeneous Einstein metrics on $S^n$ for $n=5,\ldots,9$ and
on certain products of spheres with homogeneous Einstein manifolds.
Nienhaus and Wink
\cite{NienhausWink} later used shooting to construct three non-round,
pairwise non-isometric Einstein metrics on $S^{10}$.
More recently, Buttsworth and Hodgkinson \cite{ButtsworthHodgkinson}
obtained a non-round Einstein metric on $S^{12}$ using a
computer-assisted proof.

Dancer, Hall and Wang \cite[\S5]{DancerHallWang} numerically
studied compact shrinking solitons among cohomogeneity one
metrics determined by two functions. Within the families
and parameter ranges tested, their shooting searches recovered
the Koiso--Cao soliton and known Einstein metrics.

On noncompact manifolds, Wink \cite{WinkHopf,WinkPotential}
constructed complete steady and expanding solitons.
Donovan \cite{Donovan} constructed complete non-Einstein asymptotically
conical expanding Ricci solitons in dimension four with $\SU(2)$ symmetry.
These constructions require completeness at an infinite end, whereas
our compact problem requires a second smooth collapse at finite distance.

The manifold $M$ is an $S^2$ bundle
over the six-dimensional full flag manifold $B=\SU(3)/\mathsf T^2$.
The $\SU(3)\times S^1$ action has orbit space $[0,T]$, and its principal
orbits are circle bundles over $B$. The circle fibres collapse at the two
singular orbits, each of which is diffeomorphic to $B$.
The isotropy representation of a principal orbit decomposes into
four pairwise inequivalent irreducible summands: the one-dimensional
tangent space to the circle fibre and three two-dimensional
horizontal summands. An invariant metric is therefore
determined by four functions on \([0, T]\).

For the known K\"ahler shrinker on $M$ \cite{DancerWang}, the
horizontal metric is a multiple of a fixed K\"ahler--Einstein
metric on $B$. We allow the relative scales of the horizontal
summands to vary. The proof uses a subfamily preserved by the
soliton equations in which three functions determine the metric
(see \autoref{sec:geometry}).

Smooth local solutions in this subfamily are determined by three
parameters at each singular orbit
(see \autoref{sec:singular-ivp}).
In a two-ended shooting method, the two local solutions are integrated
towards a common principal orbit. Fixing the integration distance from
one end leaves the six endpoint parameters and the integration distance
from the other to vary. The matching conditions require agreement of
the orbit metrics and their first normal derivatives, together with
the normal derivatives of the soliton potentials. These seven conditions define a
matching map whose zeros give smooth solitons on $M$
(see \autoref{sec:matching}).

We first find an approximate solution by numerical shooting.
Convergent series at the singular orbits and interval integration give
rigorous bounds on the matching map and its derivatives near the
numerical candidate. We use these bounds in a Newton--Kantorovich
argument to prove that an exact zero exists nearby.
The verification procedure and code are described in
\autoref{app:reproduction}.

Computer-assisted fixed-point arguments have been used to construct Einstein metrics
\cite{ButtsworthHodgkinson,WangNegative}, including the
K\"ahler--Einstein metric on the third del Pezzo surface
\cite{DelPezzoValidated}. For a one-parameter shooting
problem, Wang \cite{WangSpheres} used validated trajectories and
an intermediate-value argument. Cao-Labora and Rodr\'iguez-V\'azquez
\cite{CaoLaboraRodriguezVazquez} also use a finite-dimensional
matching argument in their Einstein construction.

Estimates on the potential and curvature show that the resulting
metric is non-Einstein and obstruct compatible parallel complex
structures and local product decompositions. 

\subsection*{AI disclosure}
The construction, proof arguments and computational verification
code were initially generated by ChatGPT Sol 6.1 on 30 October
2026 in response to a prompt supplied by the second-named author.
AI assistance was also used to draft and revise the manuscript
and to check the proof arguments. The authors have checked the
mathematical arguments and reviewed the verification code to
confirm that it implements the estimates used in the proof.
A complete replay of the verification program passed all
computational checks. The source code, recorded outputs and
reproduction instructions are supplied with the manuscript. The authors
take full responsibility for the manuscript and the computer-assisted proof.

\subsection*{Acknowledgements}

The second-named author gratefully acknowledges access to OpenAI's ChatGPT for
Academic Researchers program through a research workspace awarded to
Kyle Broder.  The second-named author was supported by an Australian Government
Research Training Program Scholarship.

\section{Geometric preliminaries}\label{sec:geometry}

We construct the manifold from its isotropy groups and normal discs,
following the cohomogeneity one description of Grove--Ziller
\cite[\S1]{GroveZiller} and Hoelscher \cite[\S1]{Hoelscher}.  We then
describe the invariant metrics and compute their Ricci curvature from the
homogeneous-space formula in \cite[Corollary~7.38]{Besse}.  A similar
cohomogeneity one setup was used by the first-named author to construct
non-compact gradient steady Ricci solitons in \cite{ButtsworthSteady}.

\subsection{The group diagram and the two disc bundles}

Let \(\mathsf G\) be a compact Lie group acting on a connected smooth
manifold \(M\) with cohomogeneity one.  If the orbit space is a closed
interval, the action is encoded by a group diagram
\begin{equation}\label{eq:diagram}
\begin{tikzpicture}[baseline=(current bounding box.center)]
    \node (G)  at (0,2)   {$\mathsf G$};
    \node (Km) at (-1,1)  {$\mathsf K_-$};
    \node (Kp) at (1,1)   {$\mathsf K_+$};
    \node (H)  at (0,0)   {$\mathsf H$};

    \draw[->] (Km) -- (G);
    \draw[->] (Kp) -- (G);
    \draw[->] (H) -- (Km);
    \draw[->] (H) -- (Kp);
\end{tikzpicture}
\end{equation}
or, more briefly, by
\[
  \mathsf H\subset \mathsf K_-,\mathsf K_+\subset \mathsf G.
\]
Here \(\mathsf H\) is the principal isotropy subgroup and
\(\mathsf K_\pm\) are the isotropy subgroups of the two singular orbits.
The principal and singular orbits are, respectively,
\[
  P=\mathsf G/\mathsf H,
  \qquad
  B_\pm=\mathsf G/\mathsf K_\pm.
\]
For a point on \(B_\pm\), the slice representation of \(\mathsf K_\pm\)
is transitive on the unit sphere in the normal space, with stabiliser
\(\mathsf H\).  Thus \(\mathsf K_\pm/\mathsf H\cong S^{\ell_\pm}\), where
\(\ell_\pm+1\) is the codimension of \(B_\pm\).  The slice theorem then
gives the equivariant decomposition
\begin{equation}\label{eq:disc-bundle}
  M=(\mathsf G\times_{\mathsf K_-}D_-)
    \mathop{\cup}_{\mathrm{id}_P}
    (\mathsf G\times_{\mathsf K_+}D_+),
\end{equation}
where \(D_\pm\) is the normal disc with its slice action.  Conversely,
any such group diagram determines a cohomogeneity one manifold by this
construction.

For the manifold in \autoref{thm:main}, we take
\[
  \mathsf G=\SU(3)\times S^1,
  \qquad
  \mathsf K_-=\mathsf K_+=\mathsf K=\mathsf T^2\times S^1,
\]
where \(\mathsf T^2\) is the diagonal maximal torus in \(\SU(3)\), and
\[
  \mathsf H=
  \left\{
    \left(\diag(z_1,z_2,z_3),\frac{z_1}{z_3}\right):
    z_i\in S^1,\ z_1z_2z_3=1
  \right\}.
\]
Thus \(\mathsf H\cong\mathsf T^2\), while
\(\mathsf K/\mathsf H\cong S^1\).  Both singular orbits are the full flag
manifold
\[
  \mathsf G/\mathsf K\cong \SU(3)/\mathsf T^2.
\]

\subsection{Invariant metrics and the reductive decomposition}
\label{subsec:inv-metrics}

Write
\[
  \mathfrak g=\mathfrak{su}(3)\oplus\RR Z,
\]
where \(Z\) generates the \(S^1\)-factor with period \(2\pi\).  Let
\(B_{\mathfrak{su}(3)}\) be the Killing form of
\(\mathfrak{su}(3)\).  We use the bi-invariant inner product \(Q\) on
\(\mathfrak g\) determined by
\begin{equation}\label{eq:Q-normalisation}
  Q|_{\mathfrak{su}(3)}=-B_{\mathfrak{su}(3)},
  \qquad Q(Z,Z)=1,
  \qquad Q(Z,\mathfrak{su}(3))=0.
\end{equation}
If \(E_{ij}\) denotes the elementary matrix, set
\[
  X_{ij}=\frac1{\sqrt{12}}(E_{ij}-E_{ji}),
  \qquad
  Y_{ij}=\frac{\mathrm i}{\sqrt{12}}(E_{ij}+E_{ji}),
  \qquad i<j.
\]
These six vectors are \(Q\)-orthonormal.  Together with \(Z\) and any
\(Q\)-orthonormal basis of the diagonal Cartan algebra \(\mathfrak t\),
they form a \(Q\)-orthonormal basis of \(\mathfrak g\).

Define \(\chi:\mathsf T^2\to S^1\) by
\begin{equation}\label{eq:character}
  \chi(\diag(z_1,z_2,z_3))=\frac{z_3}{z_1},
\end{equation}
and let \(\gamma:\mathfrak t\to\RR\) be determined by
\(\chi(\exp A)=\exp(\mathrm i\gamma(A))\).  Then
\[
  \mathfrak h
   =\{A-\gamma(A)Z:A\in\mathfrak t\}.
\]
We choose the following \(\Ad(\mathsf H)\)-invariant reductive
complement:
\begin{equation}\label{eq:isotropy-summands}
\begin{aligned}
  \mathfrak g&=\mathfrak h\oplus\mathfrak m,
  \qquad
  \mathfrak m=\mathfrak m_0\oplus\mathfrak m_1
                    \oplus\mathfrak m_2\oplus\mathfrak m_3,\\
  \mathfrak m_0&=\RR Z,\\
  \mathfrak m_1&=\operatorname{span}\{X_{13},Y_{13}\},\\
  \mathfrak m_2&=\operatorname{span}\{X_{12},Y_{12}\},\\
  \mathfrak m_3&=\operatorname{span}\{X_{23},Y_{23}\}.
\end{aligned}
\end{equation}
The four summands in \eqref{eq:isotropy-summands} are irreducible and pairwise
inequivalent real \(\Ad(\mathsf H)\)-modules.  The tangent space of a singular
orbit is \(\mathfrak m_1\oplus\mathfrak m_2\oplus\mathfrak m_3\), and
\(\mathfrak m_0\) is the collapsing direction. Since
\[
  [X_{ij},Y_{ij}]
     =\frac{\mathrm i}{6}(E_{ii}-E_{jj}),
\]
projection along \(\mathfrak g=\mathfrak h\oplus\mathfrak m\) gives
\begin{equation}\label{eq:vertical-brackets}
  [X_{13},Y_{13}]_{\mathfrak m}=-\frac13Z,
  \qquad
  [X_{12},Y_{12}]_{\mathfrak m}=-\frac16Z,
  \qquad
  [X_{23},Y_{23}]_{\mathfrak m}=-\frac16Z.
\end{equation}
For later use, write
\begin{equation}\label{eq:b-constants}
  (b_1,b_2,b_3)=\left(\frac13,\frac16,\frac16\right),
  \qquad [X_i,Y_i]_{\mathfrak m}=-b_iZ,
\end{equation}
where \((X_1,Y_1)=(X_{13},Y_{13})\),
\((X_2,Y_2)=(X_{12},Y_{12})\), and
\((X_3,Y_3)=(X_{23},Y_{23})\).

The following \(\Ad(\mathsf H)\)-invariant covector on \(\mathfrak m\)
determines a \(\mathsf G\)-invariant one-form \(\theta\) on
\(P=\mathsf G/\mathsf H\):
\begin{equation}\label{eq:theta-algebraic}
  \theta(Z)=1,
  \qquad
  \theta|_{\mathfrak m_1\oplus\mathfrak m_2\oplus\mathfrak m_3}=0.
\end{equation}

Schur's lemma now shows that every \(\mathsf G\)-invariant metric on
\(P\) has the form
\[
  h=f_0^2\theta^2
    +f_1^2Q|_{\mathfrak m_1}
    +f_2^2Q|_{\mathfrak m_2}
    +f_3^2Q|_{\mathfrak m_3},
  \qquad f_i>0.
\]
Now let \(g\) be a \(\mathsf G\)-invariant metric on \(M\), and choose a
geodesic meeting every principal orbit orthogonally, parametrised by
\(g\)-arc length \(t\).  Its \(\mathsf G\)-translates identify the regular
part with \(M_{\rm reg}\cong(0,T)\times P\); in this geodesic gauge, the
metric has no mixed \(\dd t\)-orbit terms, and therefore takes the form
\begin{equation}\label{eq:metric}
  g=\dd t^2+f_0(t)^2\theta^2
    +f_1(t)^2Q|_{\mathfrak m_1}
    +f_2(t)^2Q|_{\mathfrak m_2}
    +f_3(t)^2Q|_{\mathfrak m_3},
  \qquad f_i(t)>0.
\end{equation}

\subsection{Smoothness at the singular orbits}

We use the smoothness criterion of Verdiani--Ziller
\cite[Theorem~1 and \S\S3.1--3.3]{VerdianiZiller}.  In this action the
normal slice is two-dimensional, \(\mathfrak m_0\) is the only collapsing
orbit direction, and the central circle acts with weight one on the
normal slice and trivially on the tangent space of the singular orbit.
Consequently the general criterion takes a particularly simple form.

\begin{proposition}\label{prop:smooth}
Let \(s\) denote inward distance from one of the singular orbits, and
write the metric coefficients as functions of \(s\) near that orbit.
The metric in \eqref{eq:metric} and an invariant function \(u\) extend
smoothly across the orbit if and only if, for \(s\) sufficiently small,
\begin{equation}\label{eq:smoothness}
\begin{aligned}
  f_0(s)&=s\,\phi_0(s^2),
       &\phi_0(0)&=1,\\
  f_i(s)&=\phi_i(s^2),
       &\phi_i(0)&>0\quad(1\leq i\leq3),\\
  u(s)&=\phi_u(s^2),
\end{aligned}
\end{equation}
for smooth functions \(\phi_0,\phi_1,\phi_2,\phi_3,\phi_u\).
\end{proposition}

\begin{proof}
The subgroup generated by \(Z\) has period \(2\pi\), intersects
\(\mathsf H\) trivially, and generates \(\mathsf K/\mathsf H\).  Thus
the restriction of the metric to the normal disc is smooth at its centre
precisely when \(f_0\) extends as an odd function with
\(f_0'(0)=1\).  The central circle acts trivially on each tangential
summand \(\mathfrak m_i\), so the remaining metric coefficients extend
as positive even functions.  There are no mixed conditions because the
four principal-isotropy summands are pairwise inequivalent.  Finally, a
rotation-invariant smooth function on the normal disc is a smooth
function of \(s^2\).  These statements are exactly the conditions in
\eqref{eq:smoothness}.
\end{proof}

At the left-hand singular orbit these conditions apply with \(s=t\); at
the right-hand singular orbit they apply to \(f_i(T-s)\) and \(u(T-s)\),
where \(s=T-t\).  In particular,
\[
  f_0'(0)=1,
  \qquad
  f_0'(T)=-1,
  \qquad
  f_i'(0)=f_i'(T)=u'(0)=u'(T)=0
  \quad(1\leq i\leq3).
\]
We use the full parity conditions in \eqref{eq:smoothness} below.

\subsection{Ricci curvature of \texorpdfstring{$\mathsf G$}{G}-invariant metrics}

Write \(g=\dd t^2+h(t)\), where
\begin{equation}\label{eq:homogeneous-metrics}
    h(t)=f_0(t)^2\theta^2
    +f_1^2Q|_{\mathfrak m_1}
    +f_2^2Q|_{\mathfrak m_2}
    +f_3^2Q|_{\mathfrak m_3}.
\end{equation}
We first compute the Ricci curvature of the homogeneous metric \(h(t)\)
directly from the Lie brackets. For a fixed \(t\), write \(h=h(t)\),
let \(B_{\mathfrak g}\) denote the Killing form of \(\mathfrak g\), and
choose an \(h\)-orthonormal basis \(\{E_a\}\) of \(\mathfrak m\).
For $X\in\mathfrak m$, the $\Ad(\mathsf H)$-invariance of
$\mathfrak m$ gives $[X,\mathfrak h]\subset\mathfrak m$, so the
$\mathfrak h$-diagonal block of $\operatorname{ad}_X$ is zero. Hence
\[
  \tr_{\mathfrak m}(\operatorname{pr}_{\mathfrak m}\operatorname{ad}_X)
  =\tr_{\mathfrak g}(\operatorname{ad}_X)=0,
\]
where $\operatorname{pr}_{\mathfrak m}$ is projection along
$\mathfrak g=\mathfrak h\oplus\mathfrak m$, and the last equality uses
the unimodularity of the compact group $\mathsf G$.
Thus the mean-curvature vector in \cite[Corollary~7.38]{Besse}
vanishes, and the remaining formula gives, for $X\in\mathfrak m$,
\begin{equation}\label{eq:besse-homogeneous-ricci}
    \Ric_h(X,X) = \langle \mathrm{M}_\mathfrak m X, X  \rangle_h - \frac{1}{2} B_\mathfrak g(X, X),
\end{equation}
where \(\mathrm{M}_\mathfrak{m}\) is the moment map given by 
\begin{equation}\label{eq:moment-map}
\langle \mathrm{M}_\mathfrak m X, X \rangle_h
=  -\frac12\sum_a
       \bigl\|[X,E_a]_{\mathfrak m}\bigr\|_h^2+\frac14\sum_{a,b}
       \bigl\langle[E_a,E_b]_{\mathfrak m},X\bigr\rangle_h^2.
\end{equation} 

The centrality of \(Z\) gives \([Z,\mathfrak m]=0\), thus \(B_{\mathfrak g}(Z,\mathfrak g)=0\). 
Besides the brackets in \eqref{eq:vertical-brackets}, the only
remaining nonzero projected brackets are between distinct
summands $\mathfrak m_1,\mathfrak m_2,\mathfrak m_3$. For example,
\[
\begin{aligned}
  \relax[X_{13},X_{12}]&=\frac1{\sqrt{12}}X_{23},&
  [X_{13},Y_{12}]&=-\frac1{\sqrt{12}}Y_{23},\\
  [Y_{13},X_{12}]&=\frac1{\sqrt{12}}Y_{23},&
  [Y_{13},Y_{12}]&=\frac1{\sqrt{12}}X_{23}.
\end{aligned}
\]
Each of these four brackets has squared $Q$-norm $1/12$.
Cyclically permuting the matrix indices gives the same squared
norms for the other pairs of distinct summands.
Together with \eqref{eq:vertical-brackets}, these identities
provide all the bracket data needed in Besse's formula.

Since \(\Ric_h\) is \(\mathsf G\)-invariant and the four isotropy
summands are pairwise inequivalent, it is diagonal with respect to
\eqref{eq:isotropy-summands}.  Let \(R_i\) be the eigenvalue of its Ricci
endomorphism on \(\mathfrak m_i\):
\[
  \Ric_h^\sharp|_{\mathfrak m_i}=R_i\Id.
\]
\begin{proposition}
\label{prop:orbit-ricci}
The Ricci curvature of the homogeneous metric in \eqref{eq:homogeneous-metrics} is given by
\begin{equation}\label{eq:invariant-orbit-ricci}
\begin{split}
  R_0&=\frac{f_0^2}{2}\sum_{i=1}^3\frac{b_i^2}{f_i^4},\\
  R_i&=\frac1{2f_i^2}
       +\frac{f_i^4-f_j^4-f_k^4}
              {12f_1^2f_2^2f_3^2}
       -\frac{b_i^2f_0^2}{2f_i^4},
       \qquad \{i,j,k\}=\{1,2,3\},
\end{split}
\end{equation}
where \((b_1,b_2,b_3)=\left(\frac13,\frac16,\frac16\right)\).
\end{proposition}

\begin{proof}
Using the notation in \eqref{eq:b-constants}, take the
$h$-orthonormal basis of $\mathfrak m$
\[
  E_0=\frac{Z}{f_0},\qquad
  E_{i,1}=\frac{X_i}{f_i},\qquad
  E_{i,2}=\frac{Y_i}{f_i},
  \quad 1\leq i\leq3.
\]
By \eqref{eq:vertical-brackets},
\[
  [E_{i,1},E_{i,2}]_{\mathfrak m}
  =-\frac{b_i f_0}{f_i^2}E_0.
\]
Since \(Z\) is central, \(B_{\mathfrak g}(E_0,E_0)=0\) and
\([E_0,\mathfrak m]=0\). Equations \eqref{eq:besse-homogeneous-ricci}
and \eqref{eq:moment-map} therefore give
\[
  R_0
  =\frac14\sum_{a,b}
    \bigl\langle[E_a,E_b]_{\mathfrak m},E_0\bigr\rangle_h^2
  =\frac14\sum_{i=1}^3
    2\left(\frac{b_i f_0}{f_i^2}\right)^2
  =\frac{f_0^2}{2}\sum_{i=1}^3\frac{b_i^2}{f_i^4}.
\]
For \(\{i,j,k\}=\{1,2,3\}\), the same bracket identities and
the displayed brackets between distinct summands yield
\[
\begin{aligned}
  \bigl\langle
    \mathrm M_{\mathfrak m}E_{i,1},E_{i,1}
  \bigr\rangle_h
  &= -\frac{1}{2}\sum_a
  \bigl\|[E_{i,1},E_a]_{\mathfrak m}\bigr\|_h^2
  + \frac{1}{4}\sum_{a,b}
  \bigl\langle[E_a,E_b]_{\mathfrak m},E_{i,1}
  \bigr\rangle_h^2
  \\
  &=-\frac12\left(
    \frac{b_i^2f_0^2}{f_i^4}
    +\frac{f_k^2}{6f_i^2f_j^2}
    +\frac{f_j^2}{6f_i^2f_k^2}
  \right)
  +\frac{f_i^2}{12f_j^2f_k^2}\\
  &=-\frac{b_i^2f_0^2}{2f_i^4}
    +\frac{f_i^4-f_j^4-f_k^4}
           {12f_1^2f_2^2f_3^2}.
\end{aligned}
\]
Since \(B_{\mathfrak g}|_{\mathfrak m_i}=-Q|_{\mathfrak m_i}\),
\[
  -\frac12 B_{\mathfrak g}(E_{i,1},E_{i,1})
  =\frac1{2f_i^2}.
\]
Substitution into \eqref{eq:besse-homogeneous-ricci} gives
\[
\begin{aligned}
  R_i
  &=\frac1{2f_i^2}
    -\frac{b_i^2f_0^2}{2f_i^4}
    -\frac{f_k^2}{12f_i^2f_j^2}
    -\frac{f_j^2}{12f_i^2f_k^2}
    +\frac{f_i^2}{12f_j^2f_k^2}\\
  &=\frac1{2f_i^2}
    +\frac{f_i^4-f_j^4-f_k^4}{12f_1^2f_2^2f_3^2}
    -\frac{b_i^2f_0^2}{2f_i^4},
\end{aligned}
\]
as required.
\end{proof}

We now pass from the homogeneous orbit metric \(h(t)\) to the
cohomogeneity one metric
\[
g = dt^2 +h(t).
\]
Consider the normal vector field to the principal orbits \(N=\partial_t\), and define the
shape operator of the principal orbits by
\[
  LX=\nabla^g_XN,
  \qquad
  L=\frac12h^{-1}h',
\]
where \(\nabla^g\) is the Levi-Civita connection of \(g\). 
The shape operator is diagonal with \(L|_{\mathfrak m_i}=L_i\Id\) and 
\[
L_i=\frac{f_i'}{f_i}.
\]

\begin{proposition}\label{prop:ambient-ricci}
The Ricci curvature of the cohomogeneity one metric in \eqref{eq:metric} is
\begin{align*}
    \Ric_g(X,Y)
    &=\Ric_{h}(X,Y)-h(L'X,Y)
      -(\tr L)h(LX,Y),
    \\
    \Ric_g(N,X)&=0,
    \\
    \Ric_g(N,N)
    &=-(\tr L)'-\tr(L^2),
\end{align*}
for \(X, Y\) tangent to the principal orbits. 
In terms of the metric coefficients,
\[
  \Ric_g(N,N)
  =-\frac{f_0''}{f_0}
   -2\sum_{j=1}^3\frac{f_j''}{f_j},
\]
and, for unit \(X\in\mathfrak m_i\) \((0\leq i\leq3)\), 
\[
  \Ric_g(X,X)
  =
    R_i-\frac{f_i''}{f_i}
    -\frac{f_i'}{f_i}
     \left(
       \frac{f_0'}{f_0}
       +2\sum_{j=1}^3\frac{f_j'}{f_j}
       -\frac{f_i'}{f_i}
     \right),
\]
where \(R_i\) are given by \eqref{eq:invariant-orbit-ricci}.
The components between distinct isotropy summands vanish.
\end{proposition}

\begin{proof}
The normal and tangential formulas follow from
\cite[(1.2), (1.4)]{DancerWang}. The mixed components vanish
by \cite[Remark~1.20]{DancerWang}, since the isotropy
representation is multiplicity-free. Substituting
\(L_i=f_i'/f_i\) gives the formulas in metric coefficients.
\end{proof}

Since \(b_2=b_3\), the curvature formulas are symmetric under
interchange of \((f_2,L_2)\) and \((f_3,L_3)\).
Hence the Ricci equations preserve the subfamily
\(f_2=f_3\), \(L_2=L_3\). On this subfamily, 
\begin{equation}\label{eq:reduced-orbit-ricci}
\begin{split}
  R_0&=\frac{f_0^2}{18f_1^4}
       +\frac{f_0^2}{36f_2^4},\\
  R_1&=\frac1{3f_1^2}+\frac{f_1^2}{12f_2^4}
       -\frac{f_0^2}{18f_1^4},\\
  R_2=R_3&=\frac1{2f_2^2}-\frac{f_1^2}{12f_2^4}
       -\frac{f_0^2}{72f_2^4},\\
  \tr L&=L_0+2L_1+4L_2.
\end{split}
\end{equation}

\section{The singular initial value problem}\label{sec:singular-ivp}

We first derive the soliton equations on the regular set and state
the boundary value problem on the disc-bundle manifold.
We then solve the equations locally at either singular orbit.
The local existence theory for singular cohomogeneity one initial
value problems was developed by Eschenburg and Wang
\cite{EschenburgWang} in the Einstein case and extended by Buzano
\cite{Buzano} to gradient Ricci solitons. Here we use the explicit
isotropy data and smoothness conditions to obtain convergent series and
quantitative bounds for the initial solution data in the shooting argument.

\subsection{The soliton equation and boundary conditions}
We now restrict the cohomogeneity one metrics in \eqref{eq:metric} to the family $f_2=f_3$. 
Let \(u = u(t)\) be the invariant gradient soliton potential, and define the weighted mean curvature by
\begin{equation}\label{eq:shapes}
  \xi =\tr L-u'.
\end{equation}
From the Ricci curvature formulas of \autoref{sec:geometry}, the gradient soliton equation in \eqref{eq:gradient-ricci-soliton-eq} with $\lambda=1$ yields the following first-order system 
\begin{equation}\label{eq:physical}
\begin{aligned}
 f_i'&=f_iL_i,\qquad i=0,1,2,\\
 L_0'&=-\xi L_0+\frac{f_0^2}{36f_2^4}
                     +\frac{f_0^2}{18f_1^4}-1,\\
 L_1'&=-\xi L_1+\frac1{3f_1^2}+\frac{f_1^2}{12f_2^4}
                     -\frac{f_0^2}{18f_1^4}-1,\\
 L_2'&=-\xi L_2+\frac1{2f_2^2}-\frac{f_1^2}{12f_2^4}
                     -\frac{f_0^2}{72f_2^4}-1,\\
 \xi'&=-1-L_0^2-2L_1^2-4L_2^2,
\end{aligned}
\end{equation}
for \(  Y=(f_0,f_1,f_2,L_0,L_1,L_2,\xi)\).
The vector field in \eqref{eq:physical} is real analytic on
$\{f_1f_2\ne0\}$; we use the metric region $f_0,f_1,f_2>0$.

\begin{proposition}\label{prop:sufficiency}
Let $g$ be a metric of the form given in \eqref{eq:metric} on \(I\times P \subseteq M_{\rm reg}\),
with $f_2=f_3$, and let $u=u(t)$ be an invariant function.
Then $(g,u)$ satisfies \eqref{eq:gradient-ricci-soliton-eq}
with $\lambda=1$ if and only if 
\[
Y=(f_0,f_1,f_2,L_0,L_1,L_2,\xi)
\]
solves \eqref{eq:physical}.
Conversely, any solution $Y$ of \eqref{eq:physical} with
$f_0,f_1,f_2>0$ determines such a metric and a potential $u$, up to an additive constant, by
\[
  u'=L_0+2L_1+4L_2-\xi.
\]
\end{proposition}

\begin{proof}
The identities $L_i=f_i'/f_i$, for $i=0,1,2$, give the
first three equations in \eqref{eq:physical}.
By \autoref{prop:ambient-ricci}, the tangential soliton equation is equivalent to
\[
  L'=\Ric_{h(t)}^\sharp-\Id-\xi L.
\]
Using \eqref{eq:reduced-orbit-ricci}, this gives
exactly the three equations for $L_0,L_1,L_2$
in \eqref{eq:physical}.
Similarly, \autoref{prop:ambient-ricci} gives the
normal soliton equation as
\[
  -(\tr L)'-\tr(L^2)+u''=1.
\]
Since $\xi'=(\tr L)'-u''$ and $L_3=L_2$, this is equivalent to
\[
  \xi'=-1-L_0^2-2L_1^2-4L_2^2.
\]
The mixed soliton equation holds because the mixed Ricci
components vanish, as established in
\autoref{sec:geometry}, and $(\Hess_g u)(N,X)=0$.
Conversely, let $Y$ solve \eqref{eq:physical} with
$f_0,f_1,f_2>0$. Define the metric $g$ by \eqref{eq:metric}, setting
$f_3=f_2$, and choose $u$ satisfying
\[
  u'=L_0+2L_1+4L_2-\xi.
\]
The first three equations ensure that the shape operator
of $g$ has eigenvalues $L_i$, so that $\xi=\tr L-u'$.
The preceding calculations therefore show that
$\Ric_g+\Hess_g u=g$.
\end{proof}

We seek $T>0$ and a solution of \eqref{eq:physical} on $(0,T)$
with $f_0,f_1,f_2>0$ such that the metric and potential from
\autoref{prop:sufficiency} satisfy the smoothness
conditions of \autoref{prop:smooth} at both singular orbits.
\subsection{Parameters at the singular orbits and leading Taylor series}

At the left singular orbit, the smoothness conditions in
\autoref{prop:smooth} give
\[
  f_0(0)=0,\qquad f_0'(0)=1,\qquad
  f_1'(0)=f_2'(0)=u'(0)=0.
\]
We prescribe the data at this orbit
\begin{equation}\label{eq:free-data}
 \alpha_0=u''(0),\qquad
 \alpha_1=f_1(0)>0,\qquad
 \alpha_2=f_2(0)=f_3(0)>0.
\end{equation}
Fix $u(0)=0$.
Although $f_1'(0)=f_2'(0)=0$, one cannot obtain the next coefficients
by simply setting $t=0$ in \eqref{eq:physical}: both $L_0$ and $\xi$
have leading term $t^{-1}$. The finite limits needed below are
\[
  \lim_{t\to0}\frac{L_i(t)}{t}
    =\frac{f_i''(0)}{f_i(0)},
  \qquad i=1,2.
\]
They are determined by inserting the smooth even expansions of $f_1$
and $f_2$ into the equations for $L_1$ and $L_2$ and balancing their
constant terms.  For $r_1,r_2>0$, define
\begin{equation}\label{eq:lambda-values}
\begin{aligned}
 a_1(r_1,r_2)
   &=\frac1{6r_1^2}
       +\frac{r_1^2}{24r_2^4}-\frac12,\\
 a_2(r_1,r_2)
   &=\frac1{4r_2^2}
       -\frac{r_1^2}{24r_2^4}-\frac12.
\end{aligned}
\end{equation}
The resulting second derivatives are
\[
  f_i''(0)=\alpha_i a_i(\alpha_1,\alpha_2),
  \qquad i=1,2.
\]
Substituting the Taylor expansions allowed by the smoothness conditions
in \autoref{prop:smooth} into \eqref{eq:physical} and comparing
coefficients gives
\begin{equation}\label{eq:leading-series}
\begin{aligned}
 f_0(t)&=t+
   \frac{\alpha_0-2a_1-4a_2-1}{6}t^3+O(t^5),\\
 f_i(t)&=\alpha_i+\frac{\alpha_i a_i}{2}t^2+O(t^4),
 \\
 u(t)&=\frac{\alpha_0}{2}t^2+O(t^4),
\end{aligned}
\end{equation}
where $i=1,2$ and $a_i=a_i(\alpha_1,\alpha_2)$.
In \autoref{lem:local-solution}, we prove that
$(\alpha_0,\alpha_1,\alpha_2)$ uniquely determine a local
real analytic solution. Its Taylor coefficients are determined
recursively by \eqref{eq:local-recurrence}.

At the right singular orbit, we use $s=T-t$ and prescribe the data
\begin{equation}\label{eq:right-free-data}
 \omega_0=u''(T),\qquad
 \omega_1=f_1(T)>0,\qquad
 \omega_2=f_2(T)=f_3(T)>0.
\end{equation}
The functions $f_i(T-s)$ and $u(T-s)-u(T)$ have the expansions in
\eqref{eq:leading-series}, with $t$ replaced by $s$ and
$(\alpha_0,\alpha_1,\alpha_2)$ replaced by
$(\omega_0,\omega_1,\omega_2)$.
In particular,
\[
  u(T-s)=u(T)+\frac{\omega_0}{2}s^2+O(s^4).
\]

\subsection{Reduction to a regular-singular system}

The expansions in \eqref{eq:leading-series} show that $L_0$ and
$\xi$ have leading term $t^{-1}$ at the left singular orbit,
while $f_0,L_1,L_2$ vanish to at least first order.
We therefore set $x=t^2$ and write
\[
\begin{aligned}
 f_0(t)&=t\tilde f_0(x),&
 L_0(t)&=t^{-1}+t\ell_0(x),\\
 L_i(t)&=t\ell_i(x)\quad(i=1,2),&
 \xi(t)&=t^{-1}+t\hat\xi(x).
\end{aligned}
\]
Under the identities $L_i=f_i'/f_i$, $\xi=\tr L-u'$ and
$f_3=f_2$, the smoothness conditions for the metric and potential in
\autoref{prop:smooth} are equivalent to smooth extension in $x$ of
\[
  \tilde f_0,\ f_1^2,\ f_2^2,\ \ell_0,\ \ell_1,\ \ell_2,\ \hat\xi,
\]
with $\tilde f_0(0)=1$ and $f_i^2(0)>0$ for $i=1,2$.
Indeed, writing $f_0(t)=t\phi_0(x)$, $f_i(t)=\phi_i(x)$ and
$u(t)=\phi_u(x)$ as in \eqref{eq:smoothness} gives
\[
  \ell_0=2\frac{\phi_0'}{\phi_0},\qquad
  \ell_i=2\frac{\phi_i'}{\phi_i}\quad(i=1,2),\qquad
  \hat\xi=\ell_0+2\ell_1+4\ell_2-2\phi_u',
\]
where the derivatives in this display are with respect to $x$.
These expressions are smooth near $x=0$.
Conversely, positive square roots of $f_1^2$ and $f_2^2$, together
with $f_0=t\tilde f_0(t^2)$, give the required metric parity and
endpoint values. The identity defining $\xi$ gives
\[
  u'(t)=t\bigl(\ell_0+2\ell_1+4\ell_2-\hat\xi\bigr)(t^2).
\]
Consequently, for some constant $c$,
\[
  u(t)=c+\frac12\int_0^{t^2}
       \bigl(\ell_0+2\ell_1+4\ell_2-\hat\xi\bigr)(z)\,\dd z,
\]
which is a smooth function of $t^2$ and supplies the required
extension of the potential.
By \eqref{eq:physical}, the system of ODEs in these coordinates is
\begin{equation}\label{eq:regularised-system}
\begin{aligned}
 \frac{d\tilde f_0}{dx}
   &=\frac12\tilde f_0\ell_0,\\
 \frac{d(f_i^2)}{dx}
   &=f_i^2\ell_i,\qquad i=1,2,\\
 x\frac{d\ell_0}{dx}
   &=-\ell_0-\frac{\hat\xi+1}{2}
     +\frac{x}{2}\left(
       -\hat\xi\ell_0
       +\frac{\tilde f_0^2}{36f_2^4}
       +\frac{\tilde f_0^2}{18f_1^4}\right),\\
 x\frac{d\ell_1}{dx}
   &=a_1(f_1,f_2)-\ell_1
     -\frac{x}{2}\left(
       \hat\xi\ell_1
       +\frac{\tilde f_0^2}{18f_1^4}\right),\\
 x\frac{d\ell_2}{dx}
   &=a_2(f_1,f_2)-\ell_2
     -\frac{x}{2}\left(
       \hat\xi\ell_2
       +\frac{\tilde f_0^2}{72f_2^4}\right),\\
 x\frac{d\hat\xi}{dx}
   &=-\ell_0-\frac{\hat\xi+1}{2}
     -\frac{x}{2}\left(
       \ell_0^2+2\ell_1^2+4\ell_2^2\right).
\end{aligned}
\end{equation}
In \eqref{eq:regularised-system}, $a_i(f_1,f_2)$ denotes the
rational function defined in \eqref{eq:lambda-values}, applied
pointwise to the positive functions $f_1(x)$ and $f_2(x)$.
Although the new variables extend smoothly to $x=0$, the
equations remain singular there because $x$ multiplies the
derivatives of $\ell_0,\ell_1,\ell_2,\hat\xi$.

The expansions in \eqref{eq:leading-series} give
\[
 \tilde f_0(0)=1,\qquad
 f_i^2(0)=\alpha_i^2,\qquad
 \ell_i(0)=a_i(\alpha_1,\alpha_2),
 \qquad i=1,2,
\]
and
\[
 \ell_0(0)
 =
 \frac{1}{3}(\alpha_0-2a_1(\alpha_1,\alpha_2)
                   -4a_2(\alpha_1,\alpha_2)-1),
 \qquad
 \hat\xi(0)=-1-2\ell_0(0).
\]
Thus $(\alpha_0,\alpha_1,\alpha_2)$ determines the initial values of
the seven functions
\[
  (\tilde f_0,f_1^2,f_2^2,\ell_0,\ell_1,\ell_2,\hat\xi).
\]

Diagonalising the singular part of \eqref{eq:regularised-system} leads to the change of variables 
\begin{equation}\label{eq:local-coordinates}
 W=
 \begin{pmatrix}
  \tilde f_0\\
  f_1^2\\
  f_2^2\\
  \ell_0-\hat\xi\\
  \ell_1-a_1(f_1,f_2)\\
  \ell_2-a_2(f_1,f_2)\\
  2\ell_0+\hat\xi+1
 \end{pmatrix}, 
 \qquad
W_0=
 \begin{pmatrix}
  1\\
  \alpha_1^2\\
  \alpha_2^2\\
  \alpha_0-2a_1(\alpha_1,\alpha_2)-4a_2(\alpha_1,\alpha_2)\\
  0\\
  0\\
  0
 \end{pmatrix}
\end{equation}
By the smoothness conditions at the singular orbit, every component of $W$
extends smoothly to $t=0$, with $W(0)=W_0$.
The original variables are recovered from $W=(W_1,\ldots,W_7)$ by
\[
\tilde f_0=W_1, \qquad
 f_1^2 = W_2, \qquad
 f_2^2 =W_3, \qquad
 \ell_0=\frac{1}{3}(W_4+W_7-1), \qquad
 \hat\xi =\frac{1}{3}(W_7-2W_4-1)
\]
and
\[
 \ell_1=a_1(f_1,f_2)+W_5,
 \qquad
 \ell_2=a_2(f_1,f_2)+W_6.
\]
Define the rational map $\Fc$ by
\begin{equation}\label{eq:local-forcing}
 \Fc(W)=
 \begin{pmatrix}
  \dfrac{\tilde f_0\ell_0}{2}\\[6pt]
  f_1^2\ell_1\\[4pt]
  f_2^2\ell_2\\[6pt]
  -\dfrac{\hat\xi\ell_0}{2}
  +\tilde f_0^2
    \left(\dfrac1{72f_2^4}+\dfrac1{36f_1^4}\right)
  +\dfrac{\ell_0^2+2\ell_1^2+4\ell_2^2}{2}\\[10pt]
  -\dfrac{\hat\xi\ell_1}{2}
  -\dfrac{\tilde f_0^2}{36f_1^4}
  +\dfrac{\ell_1}{6f_1^2}
  +\dfrac{f_1^2(2\ell_2-\ell_1)}{24f_2^4}\\[10pt]
  -\dfrac{\hat\xi\ell_2}{2}
  -\dfrac{\tilde f_0^2}{144f_2^4}
  +\dfrac{\ell_2}{4f_2^2}
  +\dfrac{f_1^2(\ell_1-2\ell_2)}{24f_2^4}\\[10pt]
  -\hat\xi\ell_0
  +\tilde f_0^2
    \left(\dfrac1{36f_2^4}+\dfrac1{18f_1^4}\right)
  -\dfrac{\ell_0^2+2\ell_1^2+4\ell_2^2}{2}
 \end{pmatrix}.
\end{equation}
The ODE to be solved is then 
\begin{equation}\label{eq:local-system}
 x\frac{dW}{dx}
 =S(W-W_0)+x\Fc(W),
 \qquad
 S=\operatorname{diag}(0,0,0,0,-1,-1,-3/2),
\end{equation}
subject to $W(0)=W_0$, with $W_0$ as in
\eqref{eq:local-coordinates}.

\subsection{Local existence, uniqueness, and analyticity}

\begin{lemma}\label{lem:local-solution}
For each $(\alpha_0,\alpha_1,\alpha_2)\in
\RR\times(0,\infty)^2$, there exists $\varepsilon>0$ such that \eqref{eq:local-system} has a unique
smooth solution $W=W(x)$ on
$[0,\varepsilon]$ with $W(0)=W_0$.  The solution is real analytic in
$x$.  Near any fixed triple, $\varepsilon$ may be chosen
uniformly and $W$ depends real analytically, and hence smoothly, on
$(\alpha_0,\alpha_1,\alpha_2)$.
\end{lemma}

\begin{proof}
We give the fixed-point formulation because the same quantities will be
estimated explicitly below.  Put $V=W-W_0$.  Since $SW_0=0$, a regular
solution of \eqref{eq:local-system} corresponds to a fixed point $V$ of the integral operator
\begin{equation}\label{eq:volterra}
 (\mathcal TV)_j(x)=x\int_0^1 s^{-S_{jj}}
       \bigl(\Fc(W_0+V(sx))\bigr)_j\,\dd s,
  \qquad 1\le j\le7.
\end{equation}
Choose positive radii $\mathbf d=(d_1,\ldots,d_7)$ such that the box centred at $W_0$
\[
 \mathcal P(W_0,\mathbf d)
  =\{W_0+v\in\CC^7:\abs{v_j}\le d_j\}
\]
does not meet $\{W_2W_3=0\}$.  The rational map $\Fc$ is holomorphic so there are constants $B_j^{(0)}$ and $B_{jk}^{(1)}$ satisfying
\begin{equation}\label{eq:forcing-bounds}
 \abs{(\Fc)_j(W)}\le B_j^{(0)},
 \qquad
 \abs{\partial_k(\Fc)_j(W)}\le B_{jk}^{(1)}
\end{equation}
on a neighbourhood of $\mathcal{P}(W_0,\mathbf d)$.
For $\rho>0$, let $D_\rho=\{x\in\CC:|x|<\rho\}$, and let
$\mathcal A_\rho$ be the vector space of continuous maps
$V:\overline D_\rho\to\CC^7$ that are holomorphic on $D_\rho$
and satisfy $V(0)=0$, equipped with the weighted norm
\[
 \norm{V}_{\rho,\mathbf d}
  =\max_{1\le j\le7}\sup_{\abs x\le\rho}
       \frac{\abs{V_j(x)}}{d_j}.
\]
Uniform limits preserve holomorphicity on $D_\rho$ and the value at zero,
so $\mathcal A_\rho$ is a closed subspace of
$C(\overline D_\rho,\CC^7)$. The weighted norm is equivalent to the
supremum norm, hence $\mathcal A_\rho$ is a Banach space.
Since $S_{jj}\le0$, the kernels in \eqref{eq:volterra} are integrable
and
\(\int_0^1s^{-S_{jj}}\,\dd s=(1-S_{jj})^{-1}\).  Consequently, on the
closed unit ball in $\mathcal A_\rho$,
\begin{equation*}
    \norm{\mathcal TV}_{\rho,\mathbf d}\le\rho K,
 \qquad
 \|\mathcal TV-\mathcal T\widetilde V\|_{\rho,\mathbf d}
   \le\rho\Lambda
      \|V-\widetilde V\|_{\rho,\mathbf d},
\end{equation*}
where
\begin{equation}\label{eq:volterra-constants}
 K_j=\frac{B_j^{(0)}}{(1-S_{jj})d_j},
 \qquad
 \Lambda_j=
   \frac{\sum_k B_{jk}^{(1)}d_k}{(1-S_{jj})d_j},
 \qquad
 K=\max_jK_j,\qquad \Lambda=\max_j\Lambda_j.
\end{equation}
For sufficiently small $\rho$, one has $\rho K<1$ and
$\rho\Lambda<1$, so $\mathcal T$ is a strict contraction of the closed
unit ball into itself.  Its fixed point is holomorphic in $x$.  The
same choices can be made uniformly for parameters in a small complex
neighbourhood of any fixed triple, and the contraction iterates then
show that the fixed point is holomorphic as a function of
$(x,\alpha_0,\alpha_1,\alpha_2)$. For real parameters, uniqueness and complex conjugation
show that the solution is real.

Choose $0<\varepsilon<\rho$, so the analytic fixed point gives a
smooth solution on $[0,\varepsilon]$. To prove uniqueness among
smooth solutions, fix real parameters and, for $0<\delta\le\varepsilon$,
consider the real Banach space
\[
  \mathcal X_\delta
   =\{V\in C([0,\delta],\RR^7):V(0)=0\},\qquad
  \norm{V}_{\delta,\mathbf d}
   =\max_j\sup_{0\le x\le\delta}\frac{|V_j(x)|}{d_j}.
\]
The same estimates for \eqref{eq:volterra} show that $\mathcal T$
maps its closed unit ball into itself and is a contraction with
factor at most $\delta\Lambda<1$.

Let $W=W_0+V$ be any smooth solution with the prescribed initial
value. For $x>0$, \eqref{eq:local-system} gives
\[
  \frac{d}{dx}\bigl(x^{-S_{jj}}V_j(x)\bigr)
   =x^{-S_{jj}}\bigl(\Fc(W_0+V(x))\bigr)_j.
\]
Since $S_{jj}\le0$ and $V(0)=0$, the expression inside the
derivative tends to zero as $x\downarrow0$. Integrating from zero
therefore gives \eqref{eq:volterra}. After reducing $\delta$ if
necessary, both $V$ and the restriction of the analytic fixed point
lie in the closed unit ball of $\mathcal X_\delta$. Uniqueness of
the fixed point in this real function space makes them equal there.
Uniqueness for the regular ODE on $x>0$ then extends this equality
throughout their common interval of definition.
\end{proof}

For the quantitative argument we estimate the errors in truncated
Taylor series of this local solution.  Regard $W$ as a function of the
distance coordinate $t$, with $x=t^2$.
Since $SW_0=0$, \eqref{eq:local-system} gives
\begin{equation}\label{eq:regular-W}
 \frac{dW}{dt}=\mathcal G(t,W)
 :=2t\Fc(W)+\frac{2SW}{t}.
\end{equation}
Write 
\[
W=W_0+\sum_{n\ge1}W^{[n]}x^n.
\]
Comparing coefficients
in \eqref{eq:local-system} gives the recurrence
\begin{equation}\label{eq:local-recurrence}
  (n\Id-S)W^{[n]}=
  \left[\Fc\left(W_0+\sum_{j=1}^{n-1}W^{[j]}x^j\right)\right]_{n-1}.
\end{equation}
Here $W^{[n]}$ is a vector Taylor coefficient, while $W_j$ denotes
the $j$th scalar component of $W$. The notation $[\,\cdot \,]_{n-1}$
denotes the coefficient of $x^{n-1}$.
Every diagonal entry of $n\Id-S$ is positive for $n\ge1$, so this
recurrence determines $W^{[n]}$ uniquely from the preceding
coefficients. To bound the truncation error, set
\[
  P_N(x)=W_0+\sum_{n=1}^N W^{[n]}x^n.
\]
Whenever the contraction conditions above hold for $\rho$ and
$\mathbf d$, the corresponding solution satisfies
$\abs{W_j-W_{0j}}\le d_j$ on $\abs x\le\rho$.
For $0<r<\rho$, Cauchy's coefficient estimate gives
\[
  \abs{(W^{[n]})_j}
  \le r^{-n}\max_{|x|=r}|W_j(x)-W_{0j}|
  \le d_jr^{-n},\qquad n\ge1.
\]
Letting $r\uparrow\rho$ gives $\abs{(W^{[n]})_j}\le d_j\rho^{-n}$.
Summing the resulting geometric tail yields
\begin{equation}\label{eq:cauchy-tail}
 \abs{(W-P_N)_j(x)}
 \le d_j\,\frac{(\abs x/\rho)^{N+1}}{1-\abs x/\rho},
 \qquad \abs x<\rho .
\end{equation}
No holomorphic extension beyond $|x|\le\rho$ is needed.

\subsection{Numerical candidate and choice of verification data}
\label{subsec:numerical-candidate}

The initial exploration used truncated Taylor series to start the
numerical integration near each singular orbit. An adaptive
Dormand--Prince $5(4)$ solver \cite{DormandPrince} continued each
solution inward to its first zero of $\xi$.
The strict decrease of $\xi$ makes this a natural section
for comparing the two trajectories. On this section, the metric
coefficients must agree and the shape variables must have opposite
signs. A damped least-squares iteration of Levenberg--Marquardt type
\cite{Marquardt}, with a central-difference Jacobian, sought to satisfy
these matching conditions.

The initial numerical search varied the values of $f_1,f_2,f_3$
at the two singular orbits independently, without imposing
$f_2=f_3$. We checked the numerical
implementation by comparing its results for the known K\"ahler shrinking
soliton on the same bundle with an independent quadrature of that
solution \cite[\S4]{DancerWang}. Of sixteen starting points,
four approached the candidate described below, nine approached the known K\"ahler soliton, and
three reached imposed search bounds with a nonzero residual.
The new candidate exhibited $f_2=f_3$ numerically. This suggested using
the exactly invariant subfamily established in \autoref{sec:geometry}.
Details of the recorded search and of the reproducible numerical
refinement are given in \autoref{app:numerical-evidence}.

At the matching times, we compare the solution data in the coordinates
\[
 \hat Y=\mathcal C(Y)=(f_0,f_1^2,f_2^2,L_0,L_1,L_2,\xi),
\]
which avoid taking square roots when converting from $W$ by
\eqref{eq:physical-conversion}.
We parametrise the data at the two singular orbits by
\[
 \hat\alpha=(\alpha_0,\alpha_1^2,\alpha_2^2),\qquad
 \hat\omega=(\omega_0,\omega_1^2,\omega_2^2),
\]
and vary the seven parameters $p=(\hat\alpha,\hat\omega,T_R)$.
Let $p_0$ denote the exact vector used for verification.
\autoref{tab:numerical-candidate} displays its rounded coordinates
together with the fixed left comparison time $T_L$.

\begin{table}[htbp]
\centering
\caption{Numerical candidate: decimal displays of the verification
inputs, rounded to twelve decimal places.}
\label{tab:numerical-candidate}
\begin{tabular}{crcr}
\toprule
Parameter & Left value & Parameter & Right value\\
\midrule
$\alpha_0$   & $2.296157352821$ & $\omega_0$   & $-2.338419627964$\\
$\alpha_1^2$ & $0.085184676795$ & $\omega_1^2$ & $0.415303147163$\\
$\alpha_2^2$ & $0.162557980826$ & $\omega_2^2$ & $0.941298586351$\\
$T_L$       & $1.748367754142$ & $T_R$        & $1.783326986886$\\
\bottomrule
\end{tabular}
\end{table}

In particular, the choice $T_L\approx1.74836775414244$ comes from the
left trajectory's numerical $\xi=0$ hitting time. For the rigorous
calculation, we fix its stored floating-point value and interpret it as
the exact rational number
\begin{equation}\label{eq:left-time}
 T_L=\frac{1968487091515621}{2^{50}} .
\end{equation}
We allow $T_R$ to vary. Fixing $T_L$ determines
the location at which the trajectories are compared and removes
the freedom to move that location along a joined solution. The
fixed-time matching problem does not impose $\xi(T_L)=0$, although prescribing this condition and allowing both $T_L$ and $T_R$ to vary would give an alternative construction. 

The matching residual compares the left trajectory at inward time
$T_L$ with the right trajectory at inward time $T_R$. Its first
three components are the left values of $f_0,f_1^2,f_2^2$ minus
their right values; its remaining four components are the sums of
the left and right values of $L_0,L_1,L_2,\xi$.
Let $\mathcal F_{\rm num}(p)$ denote the seven-component matching
residual computed by numerically continuing the solutions from the two
singular orbits with parameters $p$. We evaluate this residual at
$p_0$ using the Taylor recurrence in \eqref{eq:local-recurrence}
and ordinary-time Taylor integration with $60$-digit decimal
arithmetic. The largest absolute component of the computed residual is
\[
 \|\mathcal F_{\rm num}(p_0)\|_\infty
 =\max_{1\le i\le7}|(\mathcal F_{\rm num}(p_0))_i|
 \approx6.57\cdot10^{-14}.
\]
\autoref{tab:numerical-convergence} shows its stability as the
integration order increases. Every row uses the same parameters
$p_0$; only the integration order changes.

Let $C$ be a fixed numerical approximation to the inverse Jacobian
of the matching residual at $p_0$; its exact rational entries are
listed in \autoref{app:preconditioner}. Starting instead with
every coordinate of $p_0$ rounded to six decimal places, the iteration
$p^{(k+1)}=p^{(k)}-C\mathcal F_{\rm num}(p^{(k)})$ gives
successive residual norms of approximately
$1.73\cdot10^{-6}$, $2.10\cdot10^{-12}$, $3.94\cdot10^{-22}$,
and $4.30\cdot10^{-31}$.
The largest absolute difference between corresponding coordinates
of the refined point and $p_0$ is approximately
$1.55\cdot10^{-13}$. The rigorous verification uses the original rational
inputs $p_0$, specified together with the matching map in
\autoref{sec:matching}.

\begin{table}[htbp]
\centering
\caption{Numerical residuals at the fixed candidate $p_0$.
The second column gives the largest absolute component of the
computed matching residual. The third gives the largest absolute
difference between corresponding components of this residual
and the reference residual.
All runs start at $t_0=1/32$ using a degree-$32$ Taylor polynomial
in $x=t^2$, and use integration steps $h\le\min(1/32,t/4)$.
The reference computation uses integration order $30$ and
$h\le\min(1/64,t/4)$.
These comparisons are numerical observations, not rigorous
error bounds.}
\label{tab:numerical-convergence}
\begin{tabular}{ccc}
\toprule
Taylor order & $\|\mathcal F_{\rm num}(p_0)\|_\infty$
 & Difference from reference\\
\midrule
$12$ & $6.57609\cdot10^{-14}$ & $4.27\cdot10^{-17}$\\
$18$ & $6.57182\cdot10^{-14}$ & $4.76\cdot10^{-25}$\\
$24$ & $6.57182\cdot10^{-14}$ & $2.32\cdot10^{-33}$\\
$30$ & $6.57182\cdot10^{-14}$ & $8.52\cdot10^{-41}$\\
\bottomrule
\end{tabular}
\end{table}

\begin{figure}[htbp]
\centering
\includegraphics[width=0.82\textwidth]{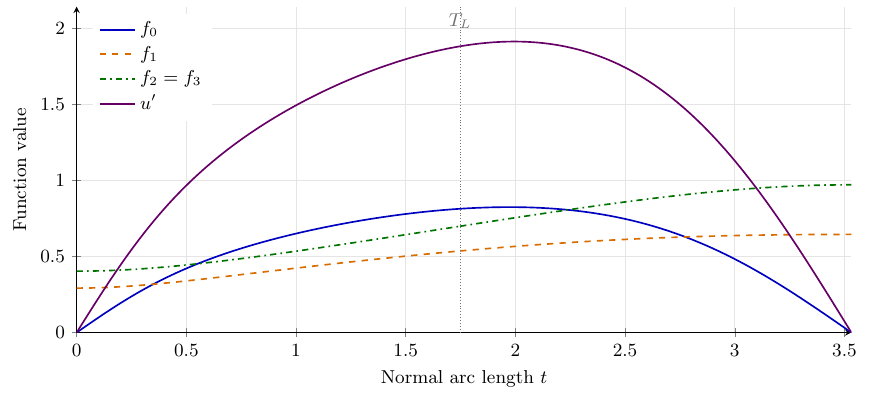}
\caption{Numerical metric functions $f_0,f_1,f_2=f_3$ and potential
derivative $u'(t)$ at the fixed candidate $p_0$, plotted on a
common scale. The soliton vector field is
$\nabla_g u=u'(t)\partial_t$, where
$u'=L_0+2L_1+4L_2-\xi$. The right trajectory, computed in inward
distance $s$, is plotted at $t=T_L+T_R-s$, with the sign of its
potential derivative reversed. The vertical grey line marks the matching
location. The two numerical pieces agree to plotting accuracy;
exact matching is proved in
\autoref{sec:matching}.}
\label{fig:numerical-candidate}
\end{figure}

\subsection{Uniform estimates and Taylor remainders}

Decimals used as input bounds denote exact rational numbers unless
explicitly identified as 64-bit floating-point values.
All root data and preconditioner entries are specified as exact rational numbers.
Using the coordinates $(\alpha_0,\alpha_1^2,\alpha_2^2)$ and
$(\omega_0,\omega_1^2,\omega_2^2)$ introduced in
\autoref{subsec:numerical-candidate} makes the initial data $W_0$
rational in these parameters.
We obtain the two chart centres by rounding the first and second triples
of the numerical candidate $p_0$ to twelve decimal places, as displayed
in \autoref{tab:numerical-candidate}:
\begin{equation}\label{eq:complex-centers}
\begin{aligned}
 \hat\alpha_{\mathrm c}
   &=(2.296157352821,0.085184676795,0.162557980826),\\
 \hat\omega_{\mathrm c}
   &=(-2.338419627964,0.415303147163,0.941298586351).
\end{aligned}
\end{equation}
For each centre in \eqref{eq:complex-centers}, allow the three
coordinates to vary independently in complex discs of
radius $\epsilon=10^{-4}$. We also use the products of three discs
of radius $\delta=10^{-6}$ about the same centres. The gap between
the radii allows Cauchy derivative estimates.
Recall the radius vector $\mathbf d=(d_1,\ldots,d_7)$ from the proof
of \autoref{lem:local-solution}: its entries specify the domain
$|W_j-W_{0j}|\le d_j$ on which the local solution is constructed. The constants $K$ and $\Lambda$
in \eqref{eq:volterra-constants} control the integral operator
and its contraction estimate. For both families of local solutions, take
\begin{equation}\label{eq:local-weights}
 \mathbf d=(1/50,1/500,1/500,1/5,1/20,1/20,1/10).
\end{equation}
These choices keep $f_1^2,f_2^2$ away from zero and permit a
common radius $\rho=1/256$ in $x$, with $\rho K<1$ and
$\rho\Lambda<1$.
Recall the bounds in \eqref{eq:forcing-bounds}: $B_j^{(0)}$
bounds the absolute value of the $j$th component of $\Fc$,
and $B_{jk}^{(1)}$ bounds the absolute value of its derivative
with respect to $W_k$.
For each family, these bounds hold uniformly on a product of
complex discs containing every $W_0+v$ with $|v_j|\le d_j$, as the
three parameters vary independently in their complex discs of radius
$\epsilon=10^{-4}$.
We refer to this product as the domain for the solution data.

\autoref{tab:local-majorants} gives rational upper bounds
for $K_j$ and $\Lambda_j$ from \eqref{eq:volterra-constants},
uniformly over the complex discs of radius $10^{-4}$.
Each bound is rounded upward to a multiple of $10^{-3}$.

\begin{table}[tb]
\centering
\caption{Bounds for the local estimates in each coordinate.
The four numerical columns give
numerators of upper bounds with common denominator $10^3$.}
\label{tab:local-majorants}
\begin{tabular}{crrrr}
\toprule
Coordinate $j$ & $K_{j,L}$ & $K_{j,R}$ & $\Lambda_{j,L}$ & $\Lambda_{j,R}$\\
\midrule
$\tilde{f}_0$ & 49388 & 20951 & 3519 & 2961\\
$f_1^2$ & 74074 & 27416 & 2827 & 10972\\
$f_2^2$ & 80726 & 143735 & 4908 & 24164\\
$\ell_0-\hat{\xi}$ & 70916 & 4884 & 7289 & 1203\\
$\ell_1-a_1(f_1,f_2)$ & 39366 & 2374 & 5407 & 341\\
$\ell_2-a_2(f_1,f_2)$ & 5080 & 610 & 1424 & 316\\
$2\ell_0+\hat{\xi}+1$ & 39956 & 2488 & 5863 & 991\\
\bottomrule
\end{tabular}
\end{table}

\begin{proposition}\label{prop:local-estimates}
On each product of three complex discs of radius $10^{-4}$,
the operator in \eqref{eq:volterra} has a unique fixed point with
$\abs{W_j-W_{0j}}\le d_j$ on $\abs{x}\le\rho=1/256$.
The constants in \eqref{eq:volterra-constants} satisfy
\[
 K_L<81,\quad K_R<144,\qquad
 \Lambda_L<8,\quad \Lambda_R<25.
\]
The real solutions have positive $f_1,f_2$ and $\tilde{f}_0\ge49/50$.
\end{proposition}
\begin{proof}
We verify the bounds needed for the contraction argument in
\autoref{lem:local-solution}, uniformly for parameters in their
complex discs of radius $10^{-4}$.
We use complex discs so that the resulting
bounds also justify the Cauchy estimates for the Taylor remainders
and their parameter derivatives.

Write
\[
 D(c,r)=\{z\in\CC:\abs{z-c}\le r\}\subset\CC,
 \qquad c\in\RR,\quad r\ge0.
\]
For products of two such discs, we have
\[
 D(c,r)\,D(c',r')
 \subseteq D\bigl(cc',\abs c\,r'+\abs{c'}r+rr'\bigr).
\]
Indeed, writing $z=c+\delta$ and $z'=c'+\delta'$ gives
$zz'-cc'=c\delta'+c'\delta+\delta\delta'$, and the triangle
inequality gives the displayed radius. Addition adds the centres
and radii. If $\abs c>r$, the disc avoids zero and its reciprocal
is exactly
\[
 \frac1{D(c,r)}
 =D\left(\frac c{c^2-r^2},\frac r{c^2-r^2}\right).
\]
All centres and radii used in the computation are rational, so
these operations can be performed in exact rational arithmetic.

For each family, we bound $W_0$ on the complex discs of
radius $10^{-4}$ by a product of discs, then increase its $j$th
coordinate radius by $d_j$ from \eqref{eq:local-weights}.
The resulting domain contains every $W_0+V(x)$ with
$\|V\|_{\rho,\mathbf d}\le1$, for all allowed data.
Since $W_2=f_1^2$ and $W_3=f_2^2$, their disc radii are
$\epsilon+d_2=\epsilon+d_3=0.0021$.
Subtracting these radii from the positive centres in
\eqref{eq:complex-centers} gives the following lower bounds
for their absolute values throughout these domains:
\[
\begin{array}{c|cc}
 &\abs{f_1^2}&\abs{f_2^2}\\ \hline
 L&0.083084676795&0.160457980826\\
 R&0.413203147163&0.939198586351 .
\end{array}
\]
The only possible poles of $\Fc$ and its derivatives occur at
$W_2W_3=0$, so these positive bounds justify the reciprocal
operations.

The routine for local estimates in \path{verify.py} applies this disc
arithmetic to the explicit rational formula in \eqref{eq:local-forcing},
using automatic differentiation to bound its first partial
derivatives as well. It returns rational bounds
\[
 \abs{(\Fc)_j(W)}\le B_j^{(0)},\qquad
 \abs{\partial_k(\Fc)_j(W)}\le B_{jk}^{(1)}
\]
valid throughout each domain. These give the
upper bounds for $K_j$ and $\Lambda_j$ recorded in
\autoref{tab:local-majorants}.
Taking the largest entry in each column and dividing by $10^3$
gives
\[
 K_L\le\frac{80726}{10^3}<81,\quad
 K_R\le\frac{143735}{10^3}<144,\quad
 \Lambda_L\le\frac{7289}{10^3}<8,\quad
 \Lambda_R\le\frac{24164}{10^3}<25.
\]
With $\rho=1/256$, the bounds on $K_L,K_R,\Lambda_L,\Lambda_R$
give, for either singular orbit,
\[
 \rho K\le\frac{144}{256}=\frac9{16}<1,\qquad
 \rho\Lambda\le\frac{25}{256}<\frac1{10}.
\]
By \eqref{eq:volterra-constants}, the first inequality ensures
that $\mathcal T$ maps the closed unit ball into itself, and the
second makes it a strict contraction. Thus each choice of data in the
discs of radius $10^{-4}$ gives a unique fixed point in this ball.
Setting $W=W_0+V$ gives
$\abs{W_j-W_{0j}}\le d_j$ on $\abs x\le1/256$.

For real parameters and real $x$, the solution is real by
\autoref{lem:local-solution}. The discs for $W_2$ and $W_3$ have
positive centres larger than their radii, so their real values
are positive; choosing the positive square roots gives $f_1,f_2>0$.
Finally, $W_{01}=1$ and $d_1=1/50$, so
\[
 \tilde f_0=W_1\ge1-\frac1{50}=\frac{49}{50}.
 \qedhere
\]
\end{proof}

For the parameter derivative estimates, restrict each parameter to the
complex disc of radius $\delta=10^{-6}$ about its centre.
For either family, write $\zeta=(\zeta_1,\zeta_2,\zeta_3)$ for
$\hat\alpha=(\alpha_0,\alpha_1^2,\alpha_2^2)$ or
$\hat\omega=(\omega_0,\omega_1^2,\omega_2^2)$.
Parameter derivatives below are taken with respect to these three
entries, in the order written.
The tail $E_N=W-P_N$ is holomorphic as a function of
$(x,\zeta_1,\zeta_2,\zeta_3)$, and \eqref{eq:cauchy-tail} holds uniformly
when the parameters lie in the discs of radius $\epsilon=10^{-4}$.
If each parameter is within $\delta$ of its centre, a circle of radius
$\eta<\epsilon-\delta$ about that parameter stays inside its disc
of radius $\epsilon$. Applying the one-variable Cauchy estimate to $E_N$, holding
the other parameters fixed, and then taking the limit as $\eta$
approaches $\epsilon-\delta=99/10^6$ from below gives
\[
  \left|\partial_{\zeta_i}(E_N)_j(x)\right|
 \le\frac{d_j}{\epsilon-\delta}
       \frac{(|x|/\rho)^{N+1}}{1-|x|/\rho},
 \qquad i=1,2,3.
\]
For the launch, choose $x_0=\rho/4$, so the geometric factor in the
tail estimate is ${4^{-N}}/{3}$. Taking $N=32$ gives the bounds for the
solution data and their parameter derivatives displayed below. Thus we use
\begin{equation}\label{eq:launch-constants}
 t_0=1/32,\qquad x_0=1/1024,\qquad N=32.
\end{equation}
Since $x_0/\rho=1/4$ and $\max_j d_j=1/5$, the common
coordinatewise bounds are
\[
 \max_j|(E_{32})_j(x_0)|\le\frac1{15\cdot2^{64}}<2^{-67},
 \qquad
  \max_{i,j}|\partial_{\zeta_i}(E_{32})_j(x_0)|
 \le\frac{10^6}{1485\cdot2^{64}}<2^{-54}.
\]
These bounds hold whenever each parameter lies in its complex
disc of radius $\delta=10^{-6}$.

To formulate and rigorously verify the matching conditions in
\autoref{sec:matching}, we use the coordinate chart $\mathcal C$
for the solution data introduced in \autoref{subsec:numerical-candidate}.
Using $f_1^2$ and $f_2^2$ as
coordinates allows conversion from $W$
without taking square roots. This conversion is the rational map
\begin{equation}\label{eq:physical-conversion}
 W\longmapsto\mathcal C(Y)=
 \bigl(t\tilde{f}_0,f_1^2,f_2^2,t^{-1}+t\ell_0,
       t\ell_1,t\ell_2,t^{-1}+t\hat{\xi}\bigr),
\end{equation}
with the substitutions from \eqref{eq:local-coordinates}.
The validated integration is carried out in the $W$ variables.
The initial interval bounds include the tails for the solution data and
their parameter derivatives given above.
At a positive terminal time, this map and its chain rule convert
the propagated interval bounds into bounds for $\mathcal C(Y)$ and its
parameter derivatives, using outward interval arithmetic.

\section{The shooting problem}\label{sec:matching}

We seek a zero of a seven-dimensional matching map $\mathcal F$
in a small box around the numerical candidate $p_0$. The proof
requires sharp interval bounds for $\mathcal F(p_0)$ and uniform
interval bounds for $D\mathcal F$ on the box. We obtain these by validated
integration of the solutions from the two singular orbits and their
parameter derivatives,
including the dependence on the right shooting time. The resulting
bounds give an exact match by a Newton--Kantorovich argument.

\subsection{The seven-dimensional matching map}

Let
\[
 Y_L(t;\alpha_0,\alpha_1,\alpha_2),\qquad
 Y_R(t;\omega_0,\omega_1,\omega_2)
\]
be the two inward local solution families supplied by
\autoref{lem:local-solution}. We use the fixed left comparison time
$T_L$ from \eqref{eq:left-time} and the rational coordinates
\begin{equation}\label{eq:parameter-chart}
 p
 =(\alpha_0,\alpha_1^2,\alpha_2^2,
     \omega_0,\omega_1^2,\omega_2^2,T_R).
\end{equation}
Its first three coordinates form $\hat\alpha$, and the next three
form $\hat\omega$. We compare the solution data through the chart
$\mathcal C$, which is injective when $f_1,f_2>0$.
Let $\mathcal D$ consist of the points $p$ with
$\alpha_1,\alpha_2,\omega_1,\omega_2,T_R>0$ for which both solutions
exist with $f_0,f_1,f_2>0$ at positive inward times through $T_L$
and $T_R$, respectively. This set is open by analytic parameter
dependence and openness of the regular ODE's existence domain.
On $\mathcal D$ define
\begin{equation}\label{eq:root-map}
 \mathcal F(p)
   =\mathcal C\bigl(Y_L(T_L;\alpha_0,\alpha_1,\alpha_2)\bigr)
    -\mathsf R\,\mathcal C
       \bigl(Y_R(T_R;\omega_0,\omega_1,\omega_2)\bigr),
 \qquad
 \mathsf R=\diag(1,1,1,-1,-1,-1,-1).
\end{equation}
The six parameters prescribed at the singular orbits and $T_R$
give seven unknowns.
We fix $T_L$ and allow $T_R$ to vary. The matching orbit
need not satisfy $\xi=0$.
Matching the shape variables and $\xi$ also matches the potential
derivatives, since $u'=L_0+2L_1+4L_2-\xi$.
Hence, we do not need a separate matching equation for $u'$.

\begin{proposition}\label{prop:closing}\label{sec:closing}
Every zero of $\mathcal F$ in $\mathcal D$ determines a smooth
$\mathsf G$-invariant gradient shrinking Ricci soliton $(g,u)$ on
the compact manifold $M$ defined in \eqref{eq:disc-bundle}.
\end{proposition}
\begin{proof}
Let $p^*$ be a zero of $\mathcal F$ and let
$T=T_L+T_R^*$.
Define the solution data on $(0,T)$ by
\[
 Y(t)=\begin{cases}
 Y_L(t;\alpha_0^*,\alpha_1^*,\alpha_2^*),&0<t\le T_L,\\
 \mathsf RY_R(T-t;\omega_0^*,\omega_1^*,\omega_2^*),
      &T_L\le t<T.
 \end{cases}
\]
Since $\mathcal F(p^*)=0$, $\mathcal C(\mathsf RY)=\mathsf R\mathcal C(Y)$,
and $\mathcal C$ is injective when $f_1,f_2>0$, the two expressions
for $Y$ agree at $t=T_L$.
Write \eqref{eq:physical} as $Y'=\mathcal V(Y)$.
Since $\mathcal V(\mathsf RY)=-\mathsf R\mathcal V(Y)$, the functions
$Y_L(t;\alpha_0^*,\alpha_1^*,\alpha_2^*)$ and
$\mathsf RY_R(T-t;\omega_0^*,\omega_1^*,\omega_2^*)$
both satisfy $Y'=\mathcal V(Y)$ as functions of $t$.
Uniqueness for this regular ODE shows that they extend to the same
smooth solution near $T_L$.

Define $u$ by integrating $u'=L_0+2L_1+4L_2-\xi$ with $u(T_L)=0$.
The potentials of the solutions constructed near the two singular orbits agree
with $u$ up to a constant. Hence $u(t)$ is smooth and even near $t=0$,
while $u(T-s)$ is smooth and even near $s=0$.
By \autoref{prop:smooth}, the metric and potential extend smoothly
across both singular orbits, so gluing the two disc bundles by
$\mathrm{id}_P$ gives a smooth $\mathsf G$-invariant pair $(g,u)$ on $M$.
The soliton equation holds on the regular set by
\autoref{prop:sufficiency} and at the singular orbits by continuity.
\end{proof}

\subsection{The box}

Using the coordinate order in \eqref{eq:parameter-chart}, we take
the box $p_0+[-r,r]^7$, with centre and radius
\begin{equation}\label{eq:root-center}
p_0=
\begin{pmatrix}
5170486699273585/2^{51}\\
6138202858770513/2^{56}\\
2928384247500027/2^{54}\\
-5265652882566845/2^{51}\\
116897443675459/2^{48}\\
4239231962734261/2^{52}\\
2007847688405285/2^{50}
\end{pmatrix},
\qquad
 r=\frac{944473296573929}{2^{73}}.
\end{equation}
Each of the seven parameters may vary by at most $r\approx10^{-7}$
from its value at $p_0$. The fractions above specify the exact values
used in the verification. \autoref{tab:numerical-candidate} gives
rounded decimal values of the entries of $p_0$.

The first three coordinates of $p_0+[-r,r]^7$ specify the allowed
data at $B_-$, and the next three specify those at $B_+$.
The computed bounds for these parameters differ from the corresponding
centre in \eqref{eq:complex-centers} by less than $2\cdot10^{-7}$
in every coordinate. Thus, the estimates near the singular orbits in
\autoref{sec:singular-ivp} apply to all these data.

To correct the numerical candidate, we use the map
$p\mapsto p-C\mathcal F(p)$, where $C$ is a fixed rational matrix
chosen to approximate $D\mathcal F(p_0)^{-1}$. Its entries are listed
in \autoref{app:preconditioner}.
We will show that, for every $p\in p_0+[-r,r]^7$, both solutions exist up to their matching times and satisfy
$f_0,f_1,f_2>0$ at positive times. This will establish
$p_0+[-r,r]^7\subset\mathcal D$, so that $\mathcal F$ is defined
throughout the box. We will also prove
\[
 \|C\mathcal F(p_0)\|_\infty<r/10,
 \qquad
 \sup_{p\in p_0+[-r,r]^7}
 \|I-CD\mathcal F(p)\|_\infty<1/10.
\]
These bounds ensure that the map sends $p_0+[-r,r]^7$
into itself. The criterion in
\autoref{subsec:contraction-criterion} then gives a unique
zero of $\mathcal F$ in $p_0+[-r,r]^7$.

\subsection{Continuation of the solutions and their derivatives}
\label{subsec:continuation}

Starting from $t_0=1/32$, we continue \eqref{eq:regular-W}
together with its first variational equations:
\[
 \partial_t W=\mathcal G(t,W),\qquad
 \partial_t(D_{\hat\alpha}W)
 =D_W\mathcal G(t,W)\,D_{\hat\alpha}W,
 \qquad \hat\alpha=(\alpha_0,\alpha_1^2,\alpha_2^2).
\]
There is no additional term because $\mathcal G$ has no explicit
parameter dependence. \autoref{sec:singular-ivp} supplies
interval bounds for $W$ and the $7\times3$ matrix
$D_{\hat\alpha}W$ at $t_0$, including the Taylor remainders and
their parameter derivatives. The same procedure applies to the
parameters $\hat\omega$ at the right singular orbit.

We use validated Taylor integration
\cite{NedialkovJacksonCorliss}, together with the standard
first-variation estimates described in
\cite[\S~2.3--2.6]{ZgliczynskiLohner}.
The computation verifies existence and positivity of
$\tilde f_0,f_1^2,f_2^2$ throughout each integration and gives
interval bounds for the solution data and their parameter
derivatives, uniformly for all allowed data at the respective singular orbit.
The solution at $p_0$ is bounded
separately to obtain a sharper bound on $\mathcal F(p_0)$.
Our implementation in \path{verify.py} and its step validation
are described in \autoref{app:validated-continuation}.

For a step from $t$ to $t+h$, let $\Phi_{t,h}(W)$ denote the
solution data at $t+h$ with initial data $W$ at $t$.
The $7\times7$ derivative $D_W\Phi_{t,h}$ propagates the
parameter derivatives by the chain rule:
\[
 D_{\hat\alpha}W(t+h,\hat\alpha)
 =D_W\Phi_{t,h}(W(t,\hat\alpha))\,
   D_{\hat\alpha}W(t,\hat\alpha).
\]
For each of the left and right families, the integration times are fixed
as the parameters vary over their prescribed ranges.
At $T_L$ and $T_{R,0}$, \eqref{eq:physical-conversion} and
the chain rule give the bounds on $\mathcal C(Y)$ and its parameter
derivatives needed for the matching map and its Jacobian.

\subsection{The variable right time}

The seventh parameter is the right shooting time. For the
right solution family, define
\[
 \hat Y_R(t;\hat\omega)
 =\mathcal C\bigl(Y_R(t;\omega_0,\omega_1,\omega_2)\bigr).
\]
In these coordinates, \eqref{eq:physical} becomes
\[
 \partial_t\hat Y_R(t;\hat\omega)
 =\hat{\mathcal V}\bigl(\hat Y_R(t;\hat\omega)\bigr).
\]
Differentiating the matching map in \eqref{eq:root-map} with respect
to $T_R$ therefore gives
\[
 \partial_{T_R}\mathcal F(p)
 =-\mathsf R\,\hat{\mathcal V}
   \bigl(\hat Y_R(T_R;\hat\omega)\bigr).
\]
To obtain interval bounds for this column and the derivatives with
respect to $\hat\omega$ for all $T_R=T_{R,0}+\sigma$, $|\sigma|\le r$,
we first continue the solutions constructed near the right singular
orbit $B_+$ to inward distance $T_{R,0}$, the seventh coordinate of $p_0$.
Let $X_R^{(0)}$ be a box containing all resulting values
$\hat Y_R(T_{R,0};\hat\omega)$ as $\hat\omega$ varies over its prescribed
ranges, and choose a larger compact box
$X_R\subset\{f_0>0,f_1^2>0,f_2^2>0\}$.
The strict inclusion
\[
 X_R^{(0)}+[-r,r]\hat{\mathcal V}(X_R)
 \subset\operatorname{int}X_R
\]
keeps every solution starting in $X_R^{(0)}$ inside $X_R$ for
both earlier and later times with $|\sigma|\le r$.
Here and below, $\hat{\mathcal V}(X_R)$ denotes a vector of
intervals computed by interval arithmetic, containing all values
of $\hat{\mathcal V}$ on $X_R$.

For initial solution data $\hat Y\in X_R^{(0)}$, let
$\hat\Psi_\sigma(\hat Y)$ be the solution data obtained by evolving this
ODE for time $\sigma$. This is the analogue of $\Phi_{t,h}$ in the
coordinates $\hat Y$ instead of $W$. No starting-time argument is
needed because the equation does not depend explicitly on time.
The computation verifies $r\Lambda_{\rm rt}<1$ and
\[
 \sup_{\hat Y\in X_R}\|D_{\hat Y}\hat{\mathcal V}(\hat Y)\|_\infty
       \le\Lambda_{\rm rt}<5.
\]
The Gronwall estimate in \eqref{eq:variational-bound}, with all
weights equal to one and in either time direction, gives
\begin{equation}\label{eq:right-time-bound}
 \norm{D_{\hat Y}\hat\Psi_\sigma(\hat Y)-I}_\infty
 \le\frac{r\Lambda_{\rm rt}}{1-r\Lambda_{\rm rt}},
 \qquad \hat Y\in X_R^{(0)},\quad \abs{\sigma}\le r.
\end{equation}
Let
\[
 P_R(\sigma,\hat\omega)
 =D_{\hat\omega}\hat Y_R(T_{R,0}+\sigma;\hat\omega).
\]
This is a $7\times3$ matrix. Bounds on $P_R(0,\hat\omega)$ are
obtained from the bounds on $D_{\hat\omega}W$ in the preceding
subsection by applying \eqref{eq:physical-conversion} and its
chain rule. The chain rule for the time extension gives
\[
 P_R(\sigma,\hat\omega)
 =D_{\hat Y}\hat\Psi_\sigma
   \bigl(\hat Y_R(T_{R,0};\hat\omega)\bigr)\,P_R(0,\hat\omega).
\]
Hence, using the induced infinity norm also for rectangular matrices,
\[
 \norm{P_R(\sigma,\hat\omega)-P_R(0,\hat\omega)}_\infty
 \le \frac{r\Lambda_{\rm rt}}{1-r\Lambda_{\rm rt}}
          \norm{P_R(0,\hat\omega)}_\infty.
\]
Let $[a_{ij},b_{ij}]$ be the computed interval containing
$(P_R(0,\hat\omega))_{ij}$ for every allowed choice of the parameters
$\hat\omega$ at $B_+$. Set
\[
 \varepsilon_R=
 \frac{r\Lambda_{\rm rt}}{1-r\Lambda_{\rm rt}}
 \max_{1\le i\le7}\sum_{j=1}^3\max\{|a_{ij}|,|b_{ij}|\}.
\]
The preceding estimate gives
\[
 a_{ij}-\varepsilon_R
 \le (P_R(\sigma,\hat\omega))_{ij}
 \le b_{ij}+\varepsilon_R,
 \qquad |\sigma|\le r.
\]
Let $\mathcal P_R$ be the $7\times3$ matrix whose $(i,j)$ entry is the interval
\[
 [a_{ij}-\varepsilon_R,\;b_{ij}+\varepsilon_R],
\]
with the lower bound rounded down and the upper bound rounded up.
Thus $\mathcal P_R$ bounds the derivative entries for every
allowed right shooting time and every allowed choice of $\hat\omega$.

Let $\mathcal P_L$ record interval bounds for the derivatives of the left
solution with respect to the parameters $\hat\alpha$ at $B_-$, evaluated
at $T_L$ and expressed in the coordinates $\hat Y$.
Thus both matrices bound derivatives in
the coordinates $\mathcal C(Y)$, whereas $\mathcal P_k$ in
\autoref{app:validated-continuation} bounds derivatives of $W$.
We collect bounds for each entry of the derivative of the matching map in
\[
 \mathcal A=
 [\,\mathcal P_L,\ -\mathsf R\mathcal P_R,
       -\mathsf R\hat{\mathcal V}(X_R)\,].
\]
Its $(i,j)$ entry is an interval containing
$\partial\mathcal F_i/\partial p_j(p)$ for every $p\in p_0+[-r,r]^7$.
The right block includes all $|\sigma|\le r$, so
$D\mathcal F(p_0+[-r,r]^7)\subset\mathcal A$.

For the later curvature estimates, we also sharpen the interval bounds
for the left solution data at $T_L$. Applying the mean-value formula
in $\hat\alpha$ over its prescribed ranges, we find that all values of the
left family at $T_L$ lie in the sum of the interval box for the solution
data at the centre and $\mathcal P_L[-r,r]^3$.

\subsection{The Newton--Kantorovich criterion}
\label{subsec:contraction-criterion}
We use the following consequence of the contraction mapping
theorem, applied to the map $p\mapsto p-CF(p)$ with $C$ fixed.

\begin{theorem}[Newton--Kantorovich criterion with a fixed inverse]
\label{lem:box-contraction}
Let $F:\mathcal U\to\RR^n$ be continuously differentiable on an open
neighbourhood of $B=p_0+[-r,r]^n$, where $r>0$.
Suppose a real $n\times n$ matrix $C$ and constants $\eta,\alpha\ge0$
satisfy
\begin{equation}\label{eq:kantorovich-hypotheses}
 \|CF(p_0)\|_\infty\le\eta,\qquad
 \sup_{p\in B}\|I-CDF(p)\|_\infty\le\alpha<1,\qquad
 \eta+\alpha r<r.
\end{equation}
Then $F$ has a unique zero $p^*$ in $B$, and
\[
 \|p^*-p_0\|_\infty\le\frac{\eta}{1-\alpha}<r.
\]
\end{theorem}
\begin{proof}
Put $\mathcal N(p)=p-CF(p)$. Integration of its derivative along
segments in the convex box gives
\[
 \|\mathcal N(p)-\mathcal N(q)\|_\infty
   \le\alpha\|p-q\|_\infty,
 \qquad
 \|\mathcal N(p)-p_0\|_\infty\le\eta+\alpha r<r.
\]
Thus $\mathcal N$ is a contraction of $B$ into its interior and
has a unique fixed point $p^*$.
Since $\|I-CDF(p_0)\|_\infty<1$, the Neumann series shows
that $CDF(p_0)$, and therefore $C$, is invertible.
Fixed points of $\mathcal N$ are consequently exactly the zeros of $F$.
Finally,
\[
 \|p^*-p_0\|_\infty
 \le\alpha\|p^*-p_0\|_\infty+\eta,
\]
which gives the stated distance bound.
\end{proof}
Let $\mathbf r_0\subset\RR^n$ be a box containing $F(p_0)$,
and let $\mathcal A\subset\RR^{n\times n}$ be a box containing
$DF(p)$ for every $p\in B$.
Recall that $\mathcal N(p)=p-CF(p)$ is the correction map.
Using interval arithmetic, define
\begin{equation}\label{eq:interval-newton}
 \mathcal E=I-C\mathcal A,\qquad
 \mathcal K=-C\mathbf r_0+\mathcal E[-r,r]^n.
\end{equation}
The entries of $\mathcal E$ bound the corresponding derivatives of
the correction map, since
\[
 D\mathcal N(p)=I-CDF(p)\in\mathcal E
 \qquad(p\in B).
\]
The box $\mathcal K\subset\RR^n$ bounds the corrected points
relative to $p_0$:
\[
 \mathcal N(B)\subset p_0+\mathcal K.
\]
Consequently, $\mathcal K\subset(-r,r)^n$ ensures that
$\mathcal N$ sends $B$ into its interior.

For an interval, its absolute value denotes the largest
absolute value it contains. We obtain bounds for the
correction at $p_0$ and the derivative throughout $B$ by taking
\[
 \eta=\max_i\abs{(-C\mathbf r_0)_i},\qquad
 \alpha=\max_i\sum_j\abs{\mathcal E_{ij}}.
\]
These satisfy
\[
 \|\mathcal N(p_0)-p_0\|_\infty\le\eta,\qquad
 \sup_{p\in B}\|D\mathcal N(p)\|_\infty\le\alpha.
\]
Finally, $0\in\mathcal E[-r,r]^n$ implies
$-C\mathbf r_0\subset\mathcal K$, and hence
\[
 \eta\le\max_i\abs{\mathcal K_i}.
\]

\subsection{Verification of the matching inequalities}

\begin{proposition}\label{prop:finite}
For every $p\in p_0+[-r,r]^7$, both solutions exist and are regular
with $f_0,f_1,f_2>0$ at positive inward times through their respective
matching times. There exist an interval box $\mathbf r_0$ containing
$\mathcal F(p_0)$ and a matrix $\mathcal A$ of intervals containing
$D\mathcal F(p_0+[-r,r]^7)$ such that, for the fixed matrix $C$,
exact rational evaluation of \eqref{eq:interval-newton} with $n=7$ gives
\begin{equation}\label{eq:terminal-bounds}
 \max_i\sum_j\abs{\mathcal E_{ij}}
 \le\alpha_*<\frac1{10},\qquad
 \max_i\frac{\abs{\mathcal K_i}}r
 \le\beta_*<\frac1{10},
\end{equation}
where the exact rational constants listed in
\autoref{app:reproduction} satisfy
\[
 \alpha_*<0.096,\qquad \beta_*<0.097.
\]
Thus the correction map $\mathcal N(p)=p-C\mathcal F(p)$ has
contraction factor less than one tenth and satisfies
\[
 \mathcal N(p_0+[-r,r]^7)\subset p_0+[-r/10,r/10]^7.
\]
At $t=T_L$, the left solution values for all allowed data at $B_-$
lie in the rational intervals in \autoref{tab:left-solution-data}.
\end{proposition}
\begin{proof}
The verification program \path{verify.py}, described in
\autoref{app:reproduction}, first checks that the prescribed parameter
ranges at both singular orbits lie within the domains used in
\autoref{sec:singular-ivp}.
\autoref{prop:local-estimates}, the Taylor recurrence in
\eqref{eq:local-recurrence} and the remainder estimates then supply
intervals containing the exact solutions and their parameter
derivatives at $t_0=1/32$. These starting bounds hold for every
allowed choice of data at either singular orbit.

At each integration step, the program verifies \eqref{eq:picard},
positivity of $\tilde f_0,f_1^2,f_2^2$, and $h\Lambda_k<7/20$.
The validated Taylor estimates in
\autoref{app:validated-continuation} guarantee that the exact
solutions exist throughout the step and give valid solution and
derivative bounds at the next step.
After $231$ left steps and
$162$ right steps, this gives bounds at $T_L$ and $T_{R,0}$,
with sharper bounds for the solutions at $p_0$.
The program checks in exact rational arithmetic that each step starts
where the preceding one ends and that the final steps end exactly at
$T_L$ and $T_{R,0}$; thus no interval of integration is skipped.

At these terminal times, \eqref{eq:physical-conversion} and the
chain rule convert the bounds on $W$ and its parameter derivatives
into bounds in the coordinates $\mathcal C(Y)$. The estimates
allowing $T_R$ to vary then cover every right shooting time in
$p_0+[-r,r]^7$. Thus both solutions exist through their
required matching times, with $f_0,f_1,f_2>0$ at positive inward
times, so
$p_0+[-r,r]^7\subset\mathcal D$.
Substituting the bounds at the matching times into the matching map in
\eqref{eq:root-map} and its derivative gives
\[
 \mathcal F(p_0)\in\mathbf r_0,\qquad
 D\mathcal F(p)\in\mathcal A
 \quad\text{for every }p\in p_0+[-r,r]^7.
\]
It remains to check that these intervals satisfy the inequalities
required by the contraction criterion.

For this final calculation, \path{verify.py} treats the lower and upper
bounds of each interval as exact fractions and uses integer numerators
and denominators. It computes $\mathcal E$ and $\mathcal K$ from
$\mathbf r_0$, $\mathcal A$ and $C$ using
\eqref{eq:interval-newton}, and verifies \eqref{eq:terminal-bounds}.
It also checks that the reported intervals contain these exact
rational results. Since $\mathcal E$ bounds $D\mathcal N$ on
$p_0+[-r,r]^7$ and $\mathcal K$ contains $\mathcal N(p_0+[-r,r]^7)-p_0$,
these inequalities prove the asserted contraction and image bounds.

Finally, the mean-value formula and the computed left
derivative bounds give intervals containing the solution values
at $t=T_L$ for all allowed data at the left singular orbit.
Rounding the lower bounds of these intervals down and their upper
bounds up gives \autoref{tab:left-solution-data}.
\end{proof}

\begin{corollary}
The matching map $\mathcal F$ has a unique zero in $p_0+[-r,r]^7$.\label{cor:root}
\end{corollary}
\begin{proof}
The initial data at $t_0$ depend analytically on the parameters.
The boxes containing the solutions throughout the integration,
including the extension in $T_R$, lie inside the regular ODE domain.
Openness of the existence domain in initial data and time therefore
makes $\mathcal F$ analytic on an open neighbourhood of $p_0+[-r,r]^7$.
By \autoref{prop:finite} and \eqref{eq:interval-newton},
\[
 \|C\mathcal F(p_0)\|_\infty\le\max_i\abs{\mathcal K_i}\le\beta_*r,
 \qquad
 \sup_{p\in p_0+[-r,r]^7}\|I-CD\mathcal F(p)\|_\infty
 \le\alpha_*.
\]
Taking $\eta=\beta_*r$ and $\alpha=\alpha_*$ gives $\alpha<1$ and
\[
 \eta+\alpha r=(\beta_*+\alpha_*)r<r/5<r.
\]
\autoref{lem:box-contraction} therefore gives the unique zero
$p^*$. Since $p^*=\mathcal N(p^*)\in p_0+\mathcal K$, we also have
\[
 \|p^*-p_0\|_\infty\le\beta_*r<r/10.
\]
\end{proof}

\begin{remark}
The uniqueness in \autoref{cor:root} holds within
$p_0+[-r,r]^7$, with the left shooting time $T_L$ fixed.
\end{remark}

The bounds at $t=T_L$ in \autoref{tab:left-solution-data} apply in
particular to the matching solution and will be used for the
curvature estimates in \autoref{sec:obstructions}.

\begin{table}[htbp]
\centering
\caption{Bounds at $t=T_L$ for all allowed data at the left singular
orbit. Divide each number in the last two columns
by $10^7$.}\label{tab:left-solution-data}
\begin{tabular}{crr}
\toprule
Coordinate&Lower numerator&Upper numerator\\
\midrule
$f_0$&8128280&8128296\\
$f_1^2$&2870430&2870434\\
$f_2^2$&4875175&4875177\\
$L_0$&1154494&1154514\\
$L_1$&2406644&2406655\\
$L_2$&3207871&3207876\\
$\xi$&-12&12\\
\bottomrule
\end{tabular}
\end{table}

\section{The manifold and the geometry of the soliton}
\label{sec:identification}

\autoref{cor:root} and \autoref{prop:closing} give
a smooth compact gradient shrinking Ricci soliton on the
disc-bundle manifold defined in \eqref{eq:disc-bundle}. We identify this manifold
and use the curvature of the soliton to prove that it is non-Einstein,
nowhere locally K\"ahler and locally irreducible.

\subsection{Identification of the manifold}

Both singular orbits are equivariantly diffeomorphic to the
six-dimensional full flag manifold
\[
  B=\mathsf G/\mathsf K\cong\SU(3)/\mathsf T^2.
\]
Define the associated complex line bundle
\begin{equation}\label{eq:line-bundle}
  \mathcal L=\SU(3)\times_\chi\CC\longrightarrow B,\qquad
  (g,w)\sim(g\vartheta,\chi(\vartheta)^{-1}w).
\end{equation}
Here $\chi:\mathsf T^2\to S^1$ is the homomorphism in
\eqref{eq:character}, which specifies how the torus rotates the complex
fibre, and $\mathcal L_{\RR}$ is the same bundle with each complex fibre
viewed as a real plane. We write $S(E)$ for the unit sphere bundle of a
real vector bundle and $\mathbb P(E)$ for the bundle of complex lines in
a complex vector bundle. The symbols $\RR$ and $\CC$ in
\eqref{eq:manifold} denote trivial line bundles over $B$.
The identifications below are fibrewise: gluing the two normal discs
gives a sphere $S^2$, which is then identified with $\CP^1$.

\begin{proposition}\label{prop:manifold-identification}
The manifold defined in \eqref{eq:disc-bundle} is equivariantly diffeomorphic to
\begin{equation}\label{eq:manifold}
  M\cong S(\RR\oplus\mathcal L_{\RR})
    \cong\mathbb P(\CC\oplus\mathcal L)
    \cong\SU(3)\times_{\mathsf T^2}\CP^1.
\end{equation}
The torus acts on $\CP^1$ by
$\vartheta\cdot[a:b]=[a:\chi(\vartheta)b]$, and the
$\mathsf G$-action in this model is
\begin{equation}\label{eq:group-action}
  (g_0,z)\cdot[g,[a:b]]=[g_0g,[a:zb]].
\end{equation}
\end{proposition}
\begin{proof}
The action of $\mathsf K$ on each normal disc is
$k\cdot w=\rho(k)w$ in a complex coordinate $w$, where
\[
  \rho:\mathsf K\longrightarrow S^1,\qquad
  \rho(\vartheta,z)=\chi(\vartheta)z,
  \qquad \ker\rho=\mathsf H.
\]
In each fibre, the boundary identification in \eqref{eq:disc-bundle}
glues two copies of the unit disc in $\CC$, each carrying the action
$w\mapsto\rho(k)w$ for $k\in\mathsf K$, along their boundary circles.
The resulting fibre is the unit sphere in $\RR\oplus\CC$, with
$\mathsf K$ acting trivially on $\RR$ and through $\rho$ on $\CC$.
Explicitly, the two maps from the unit discs are
\[
  r\eta\longmapsto
  \bigl(\pm\cos(\pi r/2),\,\sin(\pi r/2)\eta\bigr),
  \qquad 0\le r\le1,\quad |\eta|=1.
\]
They agree at $r=1$. Set $s=1-r$ on the first disc and $s=r-1$
on the second. Both maps then take the form
\[
  (s,\eta)\longmapsto
  \bigl(\sin(\pi s/2),\,\cos(\pi s/2)\eta\bigr),
\]
which is smooth across $s=0$. At $r=0$, the quotient
$\sin(\pi r/2)/r$ extends smoothly as a function of $r^2$,
with nonzero value.
Thus these maps give a smooth equivariant identification
\[
  M\cong\mathsf G\times_{\mathsf K}S(\RR\oplus\CC).
\]
To remove the extra $S^1$-factor from the group coordinate, write an
element of $\mathsf G$ as $(g,z)$ and use the map
\[
  [(g,z),(r_0,v)]\longmapsto[g,(r_0,zv)]
\]
to identify this associated sphere bundle with
$\SU(3)\times_{\mathsf T^2}S(\RR\oplus\CC)
 =S(\RR\oplus\mathcal L_{\RR})$.
The formula for $\rho$ shows that this map respects the
associated-bundle equivalence relations.

For $(r_0,v)\in S(\RR\oplus\CC)$, the image of the rank-one
Hermitian projection
\begin{equation}\label{eq:projector}
  \frac12\begin{pmatrix}
    1+r_0&\overline v\\
    v&1-r_0
  \end{pmatrix}
\end{equation}
is a complex line in $\CC^2$. Conversely, the orthogonal projection
onto any complex line has this form, giving the inverse map
$\CP^1\to S^2$. Rotation $v\mapsto zv$ corresponds to
$[a:b]\mapsto[a:zb]$, so it descends to the last two
identifications in \eqref{eq:manifold} and gives
\eqref{eq:group-action}.
\end{proof}

The singular orbits are the two sections of $M\to B$ given in each
projective fibre by $[1:0]$ and $[0:1]$. Near these sections, use the
fibre coordinates $q=b/a$ and $\widetilde q=a/b$, respectively.
An element $k\in\mathsf K$ acts by $q\mapsto\rho(k)q$ and
$\widetilde q\mapsto\rho(k)^{-1}\widetilde q$. Complex conjugation
identifies these as real representations, consistently with the two
identical real discs used in \autoref{sec:geometry}.
The principal orbit is the unit circle bundle $S(\mathcal L)$
through $(g,z)\mathsf H\mapsto[g,z]$.

\begin{proposition}\label{prop:diagram-topology}
The manifold $M$ is simply-connected. The kernel of its $\mathsf G$-action is
$\ZZ_3\times\{1\}$, so the effective group is
$\PSU(3)\times S^1$.
\end{proposition}
\begin{proof}
Since $\SU(3)$ is simply-connected and $\mathsf T^2$ is connected,
the homotopy sequence of $\mathsf T^2\to\SU(3)\to B$ gives
$\pi_1(B)=0$. The identification in \eqref{eq:manifold} exhibits $M$
as a $\CP^1$-bundle over $B$. Both base and fibre are
simply-connected, so its homotopy sequence gives $\pi_1(M)=0$.

For the action kernel, an element acting trivially on $M$ must
act trivially on $\mathsf G/\mathsf K$. The kernel of the latter
action is
\[
  \bigcap_{g\in\mathsf G}g\mathsf K g^{-1}
       =\ZZ_3\times S^1.
\]
Indeed, the intersection of all conjugates of the maximal torus in
$\SU(3)$ is its centre. On the slice, an element $(\eta\Id,z)$
of this subgroup acts as multiplication by $z$, since
$\chi(\eta\Id)=1$. Hence it acts trivially on $M$ exactly when
$z=1$, proving the assertion.
\end{proof}

\begin{remark}[The complex structure on the flag base and bundle]
We describe a complex structure on $M$ and a necessary condition for
$g$ to be K\"ahler with respect to it.
Choose the invariant complex structure $J_0$ on $B$ defined by
$J_0X_{ij}=Y_{ij}$ for $i<j$. The positive-root spaces form a
subalgebra, so $J_0$ is integrable. Let $\Omega_i$ be the invariant
two-form equal to $Q(J_0\,\cdot,\cdot)$ on $\mathfrak m_i$ and zero
whenever either argument lies in another isotropy summand.

The circle-invariant form $\theta$ from \eqref{eq:theta-algebraic}
satisfies $\theta(Z)=1$, so it is a connection one-form on $P\to B$.
With the $\Omega_i$ pulled back from $B$, \eqref{eq:vertical-brackets}
gives
\[
  \dd\theta=\tfrac13\Omega_1+\tfrac16\Omega_2+\tfrac16\Omega_3.
\]
The induced unitary connection on $\mathcal L$ has curvature
$\mathrm i\,\dd\theta$ of type $(1,1)$, so its $(0,1)$ part squares
to zero and defines a homogeneous holomorphic structure on $\mathcal L$.
The torus acts on the three complex tangent lines of $B$ by
$z_1/z_2,z_2/z_3,z_1/z_3$. Their product is $\chi^{-2}$, so the
dual determinant bundle $K_B$ has multiplier $\chi^2$, the same as
$\mathcal L^2$. Thus, with the convention in \eqref{eq:line-bundle},
\[
  K_B\cong\mathcal L^2,\qquad
  c_1(\mathcal L)=-\tfrac12c_1(B)=-[\dd\theta/(2\pi)],
\]
where $K_B$ is the canonical bundle and $c_1(B)=c_1(T^{1,0}B)$.
In particular, $\mathcal L$ is a holomorphic square root of $K_B$,
and $\mathbb P(\CC\oplus\mathcal L)$ supplies a global complex
structure on $M$.

On the regular set, lift $J_0$ to $\ker\theta$ and set
\[
  J\partial_t=f_0^{-1}Z,\qquad JZ=-f_0\partial_t.
\]
This defines an integrable complex structure for which $g$ is
Hermitian, as in \cite[\S4]{DancerWang}. Its fundamental form is
\[
  g(J\,\cdot,\cdot)
    =f_0\,\dd t\wedge\theta
      +f_1^2\Omega_1+f_2^2\Omega_2+f_3^2\Omega_3.
\]
At fixed $t$, the root brackets give
\[
  \dd_B\bigl(f_1^2\Omega_1+f_2^2\Omega_2+f_3^2\Omega_3\bigr)
       (X_{12},X_{23},Y_{13})
       =\frac{f_2^2+f_3^2-f_1^2}{\sqrt{12}},
\]
where $\dd_B$ is exterior differentiation on $B$.
For $g$ to be K\"ahler with respect to this $J$, its fundamental form
must be closed, including on each principal orbit. The preceding
calculation therefore requires $f_1^2=f_2^2+f_3^2$, or
$f_1^2=2f_2^2$ in the family $f_2=f_3$.
This condition concerns the chosen $J$; the curvature argument below
applies to all compatible parallel complex structures.
\end{remark}

\subsection{Curvature of the matching solution}\label{sec:obstructions}

Let $(g,u)$ be the soliton constructed in \autoref{prop:closing}
from the zero $p^*$ supplied by \autoref{cor:root}. We calculate its curvature
at the regular orbit $t=T_L$.
\autoref{tab:left-solution-data} gives bounds at this time for every
solution obtained from the allowed data at the left singular orbit,
including the matching solution. For the curvature calculation we enlarge the box defined by
these intervals to the following rational box:
\begin{equation}\label{eq:geometry-box}
\begin{gathered}
 f_0\in[8128,8129]/10^4,\quad
 f_1^2\in[2870,2871]/10^4,\quad
 f_2^2\in[4875,4876]/10^4,\\
 L_0\in[1154,1155]/10^4,\quad
 L_1\in[2406,2408]/10^4,\quad
 L_2\in[3207,3209]/10^4,\quad
 \xi\in[-1,1]/10^4 .
\end{gathered}
\end{equation}
Recall that
\[
 u'=L_0+2L_1+4L_2-\xi,
 \qquad
 \tr(L^2)=L_0^2+2L_1^2+4L_2^2.
\]
The orbit Ricci eigenvalues $R_0,R_1,R_2$ are given by
\eqref{eq:reduced-orbit-ricci}. We distinguish these from the Ricci
eigenvalues of the full eight-dimensional metric, denoted by
$\rho_\perp,\rho_0,\rho_1,\rho_2$ in the normal,
$\mathfrak m_0,\mathfrak m_1$, and
$\mathfrak m_2\oplus\mathfrak m_3$ directions, respectively. They are
\begin{equation}\label{eq:ambient-ricci}
\begin{aligned}
 \rho_\perp&=7-(R_0+2R_1+4R_2)-\tr(L^2)+\xi\tr L,\\
 \rho_i&=1-u'L_i,\qquad i=0,1,2,
\end{aligned}
\end{equation}
on mutually orthogonal blocks of dimensions $1,1,2,4$, respectively.
Distinct blocks may have equal eigenvalues away from the certified orbit.
To obtain the normal Ricci eigenvalue, take the trace of the
shape evolution equation:
\[
 (\tr L)'=-\xi\tr L+(R_0+2R_1+4R_2)-7.
\]
Substitution into
\[
 \rho_\perp=\Ric_g(\partial_t,\partial_t)
       =-(\tr L)'-\tr(L^2)
\]
gives the first expression in \eqref{eq:ambient-ricci}.
The tangential expressions follow from
$\Hess_g u(V,W)=u'h(t)(LV,W)$ for $V,W\in TP$.

With the convention $V\mapsto R(V,N)N$, the normal Jacobi
endomorphism is $-(L'+L^2)$. Its three tangential
block coefficients are
\begin{equation}\label{eq:jacobi}
 \kappa_i=1+\xi L_i-R_i-L_i^2,\qquad i=0,1,2.
\end{equation}
All denominators in these formulas are positive on the box in
\eqref{eq:geometry-box}.

\begin{lemma}\label{lem:curvature}
On the box in \eqref{eq:geometry-box}, the quantities in
\autoref{tab:curvature} satisfy the stated bounds.
\end{lemma}
\begin{proof}
Exact rational interval evaluation of
\eqref{eq:reduced-orbit-ricci}, \eqref{eq:ambient-ricci}, and
\eqref{eq:jacobi} on the box in \eqref{eq:geometry-box} gives the displayed bounds.
The implementation in \path{verify.py} combines repeated terms before
evaluation and checks that each interval denominator excludes zero.
The arithmetic and source are described in \autoref{app:reproduction}.
\end{proof}

\begin{table}[tb]
\centering
\caption{Exact rational curvature bounds on the box in
\eqref{eq:geometry-box}. Every row has denominator $10^6$.}
\label{tab:curvature}
\begin{tabular}{crr}
\toprule
Quantity&Lower numerator&Upper numerator\\
\midrule
$u'$&1879300&1880800\\
$\rho_\perp$&756647&760332\\
$\rho_0$&782767&783129\\
$\rho_1$&547103&547841\\
$\rho_2$&396451&397309\\
$\kappa_0$&463716&464234\\
$\kappa_1$&125155&126201\\
$\kappa_2$&10537&11042\\
\bottomrule
\end{tabular}
\end{table}

\subsection{Geometric consequences}

\begin{lemma}\label{lem:pointwise}
Let $(N_0,h)$ be a Riemannian manifold and $p\in N_0$.
Suppose a unit vector $E\in T_pN_0$ spans a simple eigenspace
of the Ricci endomorphism and the Jacobi endomorphism
\[
 \operatorname{span}\{E\}^\perp\longrightarrow \operatorname{span}\{E\}^\perp,
 \qquad X\longmapsto R(X,E)E,
\]
is invertible. No neighbourhood of $p$ admits a compatible parallel
complex structure or a nontrivial Riemannian product decomposition.
\end{lemma}
\begin{proof}
A compatible parallel complex structure commutes with the Ricci
endomorphism and hence preserves each of its eigenspaces.
Its square is $-\Id$, which is impossible on a real line.

Suppose a neighbourhood of $p$ is isometric to a nontrivial
Riemannian product. The orthogonal projection $\operatorname{pr}_1$
onto the tangent bundle of the first factor is parallel and satisfies
$0<\operatorname{rank}(\operatorname{pr}_1)<\dim N_0$.
The Ricci endomorphism is block diagonal for this product, so
$\operatorname{pr}_1$ commutes with it. Simplicity of the eigenspace spanned by
$E$ and $\operatorname{pr}_1^2=\operatorname{pr}_1$ give $\operatorname{pr}_1 E=a E$ for
$a\in\{0,1\}$.
Since $\operatorname{pr}_1$ is parallel, it also commutes with every
curvature endomorphism. Consequently, for every $X\in \operatorname{span}\{E\}^\perp$,
\[
 \operatorname{pr}_1\bigl(R(X,E)E\bigr)=R(X,E)\operatorname{pr}_1(E)=a R(X,E)E.
\]
The assumed invertibility implies that these vectors $R(X,E)E$
span $\operatorname{span}\{E\}^\perp$. Hence $\operatorname{pr}_1=a\Id$ at $p$, contradicting
the rank condition.
\end{proof}

\begin{proposition}\label{prop:all-open}
The compact soliton determined by \autoref{cor:root} has a
nonconstant potential, is non-Einstein, and has no nonempty
open K\"ahler or nontrivial product subset.
\end{proposition}
\begin{proof}
At $t=T_L$, \autoref{tab:curvature} gives
$u'(T_L)\ge1879300/10^6>0$ and
\[
  \min_{\mu\in\{\rho_0,\rho_1,\rho_2\}}
    |\rho_\perp-\mu|
    \ge\frac{22435}{10^6}>0.
\]
Thus the potential is nonconstant and the metric is non-Einstein.
Moreover, for every vector $V$ tangent to this orbit,
\[
  g(R(V,\partial_t)\partial_t,V)
    \ge\frac{10537}{10^6}\,g(V,V).
\]
Thus the normal Ricci eigenspace is simple and the normal
Jacobi endomorphism is positive definite.

On the connected regular interval $(0,T)$, the solution of
\eqref{eq:physical} is analytic. The six functions
\[
 \rho_\perp-\rho_0,\quad \rho_\perp-\rho_1,\quad
 \rho_\perp-\rho_2,\quad \kappa_0,\quad \kappa_1,\quad \kappa_2
\]
are analytic there. Their product is nonzero at $T_L$ by
\autoref{tab:curvature}, so its zero set is discrete in $(0,T)$.
At every time outside that set, $\partial_t$ spans a simple
Ricci eigenspace, and the normal Jacobi endomorphism is
invertible: its three block values are $\kappa_0,\kappa_1,\kappa_2$.
\autoref{lem:pointwise} applies at every point of such an orbit.

Let $U\subset M$ be any nonempty open subset.
It meets the regular set, since the two singular orbits have
codimension two. Arc length is a submersion on the regular
set, so the image of this intersection contains an open
interval. That interval contains a time outside the discrete
exceptional set, and $U$ therefore contains a point
to which \autoref{lem:pointwise} applies. A K\"ahler structure
or product decomposition on $U$ would restrict to a
neighbourhood of that point, which is impossible.
\end{proof}

\begin{proof}[Proof of \autoref{thm:main}]
\autoref{cor:root} supplies the matching solution, and
\autoref{prop:closing} constructs the smooth compact metric
and potential on the prescribed double disc bundle.
\autoref{prop:diagram-topology} proves that $M$ is simply-connected.
\autoref{prop:manifold-identification} identifies the total
space and its two singular orbits. \autoref{prop:all-open}
proves that the metric is non-Einstein, nowhere locally K\"ahler
and locally irreducible.
\end{proof}

\appendix
\section{Reproducibility and the finite certificate}
\label{app:reproduction}

\subsection{Arithmetic and implementation}

The recorded replay used Python~3.12.10, NumPy~2.3.5 and Node~22.23.4
on 64-bit AMD64 Windows. Interval bounds use binary64 arithmetic
with round-to-nearest, gradual underflow, and correctly rounded
addition, subtraction, multiplication and division. Each elementary
result is enlarged by one \texttt{nextafter} step towards each relevant
infinity. Products use all four endpoint products, reciprocals use
$[1/b,1/a]$ when $0\notin[a,b]$, and matrix products use explicit
interval sums. Rational coefficients are enclosed by outward-rounded
division of exactly represented integers; the root data and
preconditioner entries are the exact dyadic rationals printed in the manuscript.

Taylor coefficients are computed by standard power-series arithmetic,
and derivatives by the product and quotient rules applied to the
underlying functions. The remainder estimates in
\autoref{sec:singular-ivp} and \autoref{app:validated-continuation}
give bounds for the solutions and their parameter derivatives.
Every proposed integration step must pass the interval tests in
\autoref{app:validated-continuation}; overflow, invalid operations
and division by zero cause rejection.

For the final checks, \texttt{Fraction.from\_float} converts interval
endpoints to exact dyadic rationals. The checker recomputes the matrix
products and fixed-point image using rational arithmetic, checks the
stated containments and inequalities, and verifies that consecutive
integration steps meet exactly. Its constants in \autoref{prop:finite} are
\begin{align}
 \alpha_*&=
 \frac{247976859883368044686315079332527}
      {2^{111}},\label{eq:alpha}\\
 \beta_*&=
 \frac{236696180595019440104046978422827173381464177735}
      {944473296573929\cdot2^{111}}.
 \label{eq:beta}
\end{align}
They satisfy $\alpha_*<0.096$ and $\beta_*<0.097$.

\subsection{Validated continuation}
\label{app:validated-continuation}

We record the implementation of the standard validated Taylor
estimates used in \autoref{subsec:continuation}; see
\cite{NedialkovJacksonCorliss} and
\cite[Lemma~1 and Sections~2.4--2.6]{ZgliczynskiLohner}.
The nonautonomous equation is covered by adjoining the time
variable with derivative $1$. Below, derivatives with respect to
$W$ keep the initial time fixed.

At $t=t_k$, let $c_k$ approximate the solution at $p_0$.
The vectors $r_k$ and $r_k^{(0)}$ bound, componentwise,
the differences from $c_k$ for all allowed data at the relevant singular
orbit and for the data at $p_0$, respectively. The $7\times3$
matrix $\mathcal P_k$ records an interval bound for each parameter derivative.
For a step of length $h=t_{k+1}-t_k$, the program finds a compact
box $X_k$, with $\tilde f_0,f_1^2,f_2^2>0$, and checks
\begin{equation}\label{eq:picard}
 c_k+[-r_k,r_k]+[0,h]\mathcal G([t,t+h],X_k)
 \subset\operatorname{int}X_k.
\end{equation}
It uses the weighted maximum norm and its induced matrix norm,
\begin{equation}\label{eq:flow-weights}
 w=(1,1/10,1/10,10,1,1,10),\qquad
 \norm{v}_w=\max_i\frac{\abs{v_i}}{w_i},\qquad
 \norm{A}_w=\max_i\sum_j\abs{A_{ij}}\frac{w_j}{w_i},
\end{equation}
and bounds $\norm{D_W\mathcal G(\tau,W)}_w\le\Lambda_k$
for $\tau\in[t,t+h]$ and $W\in X_k$, requiring
$h\Lambda_k<7/20$. Since $t>0$ and the first three coordinates
are positive, $\mathcal G$ is analytic near this box.
The strict inclusion in \eqref{eq:picard} guarantees continuation
inside $X_k$. Gronwall's inequality for the first variational
equation gives, for $W_*\in c_k+[-r_k,r_k]$ and $0\le\sigma\le h$,
\begin{equation}\label{eq:variational-bound}
 \norm{D_W\Phi_{t,\sigma}(W_*)-I}_w
 \le e^{h\Lambda_k}-1
 \le\frac{h\Lambda_k}{1-h\Lambda_k}.
\end{equation}
Thus the matrix whose $(i,j)$ entry is
\[
 \boldsymbol\Sigma_{ij}
 =I_{ij}+\frac{h\Lambda_k}{1-h\Lambda_k}\frac{w_i}{w_j}[-1,1]
\]
bounds the flow derivative throughout the step.

Let $T_n(t,W)$ be the coefficient of $v^n$ in
$\Phi_{t,v}(W)$. The program computes these coefficients and
their derivatives by interval evaluation of the recurrence
\[
 T_0(t,W)=W,\qquad
 T_{n+1}(t,W)=\frac1{n+1}[v^n]\,
 \mathcal G\left(t+v,\sum_{j=0}^{n}T_j(t,W)v^j\right).
\]
Here $[v^n]$ extracts the coefficient of $v^n$ in each coordinate.
The finite Taylor estimates of \cite[Lemma~1]{ZgliczynskiLohner} give
\[
\begin{aligned}
 \Phi_{t,h}(c_k)&\in
 \sum_{n=0}^{14}h^nT_n(t,c_k)
       +h^{15}T_{15}([t,t+h],X_k),\\
 \mathcal J_k&=
 \sum_{n=0}^{8}h^nD_WT_n(t,c_k+[-r_k,r_k])
       +h^9\mathcal R_9,
\end{aligned}
\]
where $\mathcal J_k$ bounds each entry of
$D_W\Phi_{t,h}(c_k+[-r_k,r_k])$.
The remainder $\mathcal R_9$ bounds
$D_WT_9(\tau,W)\boldsymbol\Sigma$ for
$\tau\in[t,t+h]$ and $W\in X_k$. This bound is computed by applying
the product and quotient rules to the same recurrence, with the
initial perturbation matrix at $T_0=W$ set to $\boldsymbol\Sigma$.

Choose $c_{k+1}$ as the midpoint of the computed interval for
$\Phi_{t,h}(c_k)$, and let $q_k$ bound its componentwise radius
about $c_{k+1}$, including rounding error. The mean-value formula
and the chain rule justify the updates
\[
 r_{k+1}\ge q_k+\abs{\mathcal J_k}r_k,\qquad
 r_{k+1}^{(0)}\ge q_k+\abs{\mathcal J_k}r_k^{(0)},\qquad
 \mathcal P_{k+1}\supset\mathcal J_k\mathcal P_k.
\]
Here $\abs{\mathcal J_k}$ denotes the componentwise supremum of
the absolute values of its entries. All bounds are rounded
outwards, and $0\le r_k^{(0)}\le r_k$ is checked at every step.
The mesh is common to all allowed data at a given
singular orbit and is held fixed when differentiating. Stored times and step sizes
are interpreted as exact dyadic rationals, and each equality
$t_k+h=t_{k+1}$ is checked in exact fraction arithmetic, including
the final step to the prescribed matching time.

\subsection{Numerical exploration and reproduction}
\label{app:numerical-evidence}

The historical search described in \autoref{subsec:numerical-candidate}
is recorded in the embedded \path{numerical-discovery.json}, which
contains the initial and final parameters, outcomes and source fingerprints.
The historical record is not an input to the rigorous matching test.

The optional command
\begin{lstlisting}
python -B verify.py --numerics --output numerical-evidence
\end{lstlisting}
reproduces the refinement, convergence tables and figure data using
$60$-digit decimal arithmetic. Starting from the six-decimal rounding
of $p_0$, the fixed-inverse iteration uses order $24$ and a maximum
step size of $1/32$, stopping when the residual falls below $10^{-28}$.
A separate order-$30$ evaluation with a maximum step size of $1/64$ gives a 
residual of approximately $4.32\cdot10^{-31}$.
The convergence table uses the same exact verification inputs at every
order. This mode requires Python and NumPy and reports numerical
evidence only; full interval verification is described below.

\subsection{Reproduction of the computation}

The sole computational attachment, \path{verify.py}, contains all numerical inputs,
readable source blocks and the historical search record, and checks
each block's hash before extraction. Its SHA256, using UTF-8 with
LF line endings, is
\begin{quote}\small
\path{4D9C89E9945DE6FEC40C268EDC9EA344AE4A1909115DA68AD3F5231C234BFB65}.
\end{quote}
The source uses documented ASCII versions of the manuscript's notation.
States, radii, flow weights and flow Jacobians follow $(W_1,\ldots,W_7)$;
solution data follow $\mathcal C(Y)$, and parameters follow
\eqref{eq:parameter-chart}. The tuple \texttt{p\_hat} is
$\hat\alpha$ or $\hat\omega$, denoted by $\zeta$ in the derivative estimates.
The Nomizu calculation uses tangent-block order
$(\mathfrak m_0,\mathfrak m_1,\mathfrak m_2,\mathfrak m_3)$ and
polynomial-variable order $(f_0^2,f_1^2,f_2^2,f_3^2)$.
The source blocks have the following roles:
\begin{center}
\begin{tabular}{@{}lp{0.64\linewidth}@{}}
\toprule
Block & Mathematical role\\
\midrule
\texttt{CAP} & Complex-disc majorants\\
\texttt{FLOW} & Local launch, validated flow and parameter derivatives\\
\texttt{ROOT\_TEST}, \texttt{EXACT} & Interval matching and exact terminal tests\\
\texttt{NOMIZU}, \texttt{GEOMETRY} & Orbit Ricci formulas and curvature bounds\\
\texttt{NUMERICS} & Numerical refinement and convergence study\\
\bottomrule
\end{tabular}
\end{center}

With the runtime described above and Node on the executable search path,
run
\begin{lstlisting}
python -B verify.py --output replay
\end{lstlisting}
The directory must be new or empty; it receives the extracted sources,
step traces, certificates and a summary identifying the executed
verifier by its SHA256. All flow bounds are recomputed, and a successful
full run verifies \eqref{eq:alpha}--\eqref{eq:beta} and reports
\texttt{FULL COMPUTATIONAL PROOF REPLAY PASS}.
The wrapper rejects Python optimisation (\texttt{-O}, \texttt{-OO}),
since assertions perform mathematical checks.
The optional \texttt{--quick} mode checks source integrity, local
estimates and exact curvature, but omits the integrations and does
not verify the matching zero. Omitting \texttt{--output} uses a
temporary directory; \texttt{--node} specifies the Node executable by path.

\section{The exact preconditioner}\label{app:preconditioner}

The table specifies every entry of the fixed $7\times7$
matrix $C$ in \autoref{prop:finite}.
The column order is the order of the solution data
$(f_0,f_1^2,f_2^2,L_0,L_1,L_2,\xi)$ and the row order is
the coordinate order of \eqref{eq:parameter-chart}.
All entries are exact dyadic rationals.
The matrix was obtained by numerically inverting a central-difference
Jacobian at the candidate. Its entries are then treated as fixed exact
inputs. Its adequacy is established by the bound on $I-C\mathcal A$
in \eqref{eq:terminal-bounds}.
{
\renewcommand{\arraystretch}{1.6}
\setlength{\tabcolsep}{5pt}
\begin{longtable}{@{}rr >{$\displaystyle}l<{$}
                 @{\hspace{1cm}}
                 rr >{$\displaystyle}l<{$}
                 @{\hspace{1cm}}
                 rr >{$\displaystyle}l<{$}@{}}
\caption{Exact entries of $C$.}\label{tab:preconditioner}\\
\toprule
$i$ & $j$ & \multicolumn{1}{c}{$C_{ij}$}
& $i$ & $j$ & \multicolumn{1}{c}{$C_{ij}$}
& $i$ & $j$ & \multicolumn{1}{c}{$C_{ij}$}\\
\midrule
\endfirsthead
\toprule
$i$ & $j$ & \multicolumn{1}{c}{$C_{ij}$}
& $i$ & $j$ & \multicolumn{1}{c}{$C_{ij}$}
& $i$ & $j$ & \multicolumn{1}{c}{$C_{ij}$}\\
\midrule
\endhead
\bottomrule
\endfoot
1&1&-\frac{1662985824653601}{2^{53}}
&3&4&\frac{3719421233646797}{2^{56}}
&5&6&\frac{5206164865307723}{2^{54}}\\[4pt]
1&2&-\frac{5848244028560447}{2^{50}}
&3&5&-\frac{6714261369249611}{2^{56}}
&5&7&-\frac{3729814050783367}{2^{54}}\\[4pt]
1&3&-\frac{6928056599632281}{2^{51}}
&3&6&-\frac{8846646475250993}{2^{56}}
&6&1&\frac{4885334077849001}{2^{54}}\\[4pt]
1&4&\frac{1161097438931889}{2^{58}}
&3&7&\frac{8620254260899073}{2^{57}}
&6&2&\frac{5911460039915759}{2^{54}}\\[4pt]
1&5&\frac{1707562811852043}{2^{49}}
&4&1&-\frac{2158710640032255}{2^{50}}
&6&3&-\frac{10662693769297}{2^{42}}\\[4pt]
1&6&\frac{5294222100703669}{2^{51}}
&4&2&\frac{2757281149007953}{2^{50}}
&6&4&-\frac{7681328065984249}{2^{55}}\\[4pt]
1&7&-\frac{7375255286338991}{2^{52}}
&4&3&\frac{4748025778543905}{2^{50}}
&6&5&\frac{2088501715118339}{2^{52}}\\[4pt]
2&1&\frac{6478423952897643}{2^{57}}
&4&4&\frac{6469431493319919}{2^{52}}
&6&6&-\frac{7706576135111947}{2^{53}}\\[4pt]
2&2&\frac{5244510445574699}{2^{53}}
&4&5&-\frac{4846042474695305}{2^{52}}
&6&7&-\frac{6706965671855907}{2^{54}}\\[4pt]
2&3&-\frac{191553720642767}{2^{49}}
&4&6&\frac{692484341901631}{2^{50}}
&7&1&\frac{3655858895168945}{2^{53}}\\[4pt]
2&4&-\frac{5077945051062151}{2^{57}}
&4&7&\frac{5374903453183357}{2^{53}}
&7&2&\frac{5342405482571519}{2^{52}}\\[4pt]
2&5&-\frac{1386291986779115}{2^{52}}
&5&1&-\frac{6577677915096935}{2^{55}}
&7&3&\frac{6210195079308597}{2^{51}}\\[4pt]
2&6&\frac{3981448239452333}{2^{54}}
&5&2&-\frac{5711355631945835}{2^{51}}
&7&4&-\frac{3084849707999931}{2^{54}}\\[4pt]
2&7&\frac{4946699235090279}{2^{58}}
&5&3&-\frac{3261102353781187}{2^{53}}
&7&5&\frac{7102307323814831}{2^{54}}\\[4pt]
3&1&\frac{2968208405014335}{2^{55}}
&5&4&-\frac{5195675803044451}{2^{54}}
&7&6&-\frac{2567705460612311}{2^{51}}\\[4pt]
3&2&\frac{4515551837545793}{2^{55}}
&5&5&-\frac{5590808085627475}{2^{54}}
&7&7&\frac{2379749543428999}{2^{54}}\\[4pt]
3&3&\frac{1747820534410151}{2^{53}}
& & & 
& & & \\[4pt]
\end{longtable}
}


\begin{thebibliography}{BHH+26}

\bibitem[ALB87]{Besse}
A.~L. Besse,
\emph{Einstein Manifolds},
Ergebnisse der Mathematik und ihrer Grenzgebiete (3), vol.~10,
Springer-Verlag, Berlin, 1987.

\bibitem[Boh98]{Bohm}
C.~B\"ohm,
\newblock Inhomogeneous Einstein metrics on low-dimensional spheres and
other low-dimensional spaces,
\newblock \emph{Invent. Math.} \textbf{134} (1998), no.~1, 145--176.
\newblock \href{https://doi.org/10.1007/s002220050261}
{doi:10.1007/s002220050261}.

\bibitem[But25]{ButtsworthSteady}
T.~Buttsworth,
\newblock $\SU(2)$-invariant steady gradient Ricci solitons on four-manifolds,
\newblock \emph{Comm. Anal. Geom.} \textbf{33} (2025), no.~2, 365--401.
\newblock \href{https://doi.org/10.4310/CAG.250530011944}
{doi:10.4310/CAG.250530011944}.

\bibitem[BHH+26]{DelPezzoValidated}
T.~Buttsworth, W.~Hadden, E.~Heyes, D.~Platt, and T.~Wiseman,
\newblock The K\"ahler--Einstein metric on the third del Pezzo surface,
\newblock preprint, \href{https://arxiv.org/abs/2609.32655}
{arXiv:2609.32655} (2026).

\bibitem[BH26]{ButtsworthHodgkinson}
T.~Buttsworth and L.~Hodgkinson,
\newblock Computationally assisted proof of a novel
$\mathsf O(3)\times\mathsf O(10)$-invariant Einstein metric on $S^{12}$,
\newblock \emph{J. Lond. Math. Soc.} \textbf{113} (2026), no.~2,
Article e70477.

\bibitem[Buz11]{Buzano}
M.~Buzano,
\newblock Initial value problem for cohomogeneity one gradient Ricci
solitons,
\newblock \emph{J. Geom. Phys.} \textbf{61} (2011), no.~6, 1033--1044.

\bibitem[Cao96]{Cao}
H.-D.~Cao,
\newblock Existence of gradient K\"ahler--Ricci solitons,
\newblock in \emph{Elliptic and Parabolic Methods in Geometry},
A K Peters, Wellesley, MA, 1996, pp.~1--16.

\bibitem[CLRV26]{CaoLaboraRodriguezVazquez}
G.~Cao-Labora and A.~Rodr\'iguez-V\'azquez,
\newblock Inhomogeneous Einstein metrics on complex projective spaces,
\newblock preprint, \href{https://arxiv.org/abs/2608.16880}
{arXiv:2608.16880} (2026).

\bibitem[DHW13]{DancerHallWang}
A.~S.~Dancer, S.~J.~Hall, and M.~Y.~Wang,
\newblock Cohomogeneity one shrinking Ricci solitons: an analytic and
numerical study,
\newblock \emph{Asian J. Math.} \textbf{17} (2013), no.~1, 33--62.

\bibitem[DW11]{DancerWang}
A.~S.~Dancer and M.~Y.~Wang,
\newblock On Ricci solitons of cohomogeneity one,
\newblock \emph{Ann. Global Anal. Geom.} \textbf{39} (2011), no.~3,
259--292.

\bibitem[Der12]{Derdzinski}
A.~Derdzinski,
\newblock Ricci solitons,
\newblock \emph{Wiadomo\'sci Matematyczne} \textbf{48} (2012), no.~1,
1--32.

\bibitem[Don25]{Donovan}
P.~Donovan,
\newblock Cohomogeneity one 4-dimensional gradient Ricci solitons,
\newblock \emph{J. Geom. Anal.} \textbf{35} (2025), Article 401.
\newblock \href{https://doi.org/10.1007/s12220-025-02220-z}
{doi:10.1007/s12220-025-02220-z}.

\bibitem[DP80]{DormandPrince}
J.~R.~Dormand and P.~J.~Prince,
\newblock A family of embedded Runge--Kutta formulae,
\newblock \emph{J. Comput. Appl. Math.} \textbf{6} (1980), no.~1, 19--26.
\newblock \href{https://doi.org/10.1016/0771-050X(80)90013-3}
{doi:10.1016/0771-050X(80)90013-3}.

\bibitem[EW00]{EschenburgWang}
J.-H.~Eschenburg and M.~Y.~Wang,
\newblock The initial value problem for cohomogeneity one Einstein metrics,
\newblock \emph{J. Geom. Anal.} \textbf{10} (2000), no.~1, 109--137.
\newblock \href{https://doi.org/10.1007/BF02921808}
{doi:10.1007/BF02921808}.

\bibitem[GZ00]{GroveZiller}
K.~Grove and W.~Ziller,
\newblock Curvature and symmetry of Milnor spheres,
\newblock \emph{Ann. of Math.} (2) \textbf{152} (2000), no.~1, 331--367.

\bibitem[Hoe10]{Hoelscher}
C.~A.~Hoelscher,
\newblock Classification of cohomogeneity one manifolds in low dimensions,
\newblock \emph{Pacific J. Math.} \textbf{246} (2010), no.~1, 129--185.

\bibitem[Koi90]{Koiso}
N.~Koiso,
\newblock On rotationally symmetric Hamilton's equation for
K\"ahler--Einstein metrics,
\newblock in \emph{Recent Topics in Differential and Analytic Geometry},
Adv. Stud. Pure Math., vol.~18-I, Academic Press, Boston, 1990,
pp.~327--337.

\bibitem[Mar63]{Marquardt}
D.~W.~Marquardt,
\newblock An algorithm for least-squares estimation of nonlinear parameters,
\newblock \emph{J. Soc. Indust. Appl. Math.} \textbf{11} (1963), no.~2,
431--441.
\newblock \href{https://doi.org/10.1137/0111030}
{doi:10.1137/0111030}.

\bibitem[NJC99]{NedialkovJacksonCorliss}
N.~S.~Nedialkov, K.~R.~Jackson, and G.~F.~Corliss,
\newblock Validated solutions of initial value problems for ordinary
differential equations,
\newblock \emph{Appl. Math. Comput.} \textbf{105} (1999), no.~1, 21--68.

\bibitem[NW25]{NienhausWink}
J.~Nienhaus and M.~Wink,
\newblock Einstein metrics on the ten-sphere,
\newblock \emph{J. Eur. Math. Soc.} (2025), published online first.
\newblock \href{https://doi.org/10.4171/JEMS/1668}
{doi:10.4171/JEMS/1668}.

\bibitem[Per02]{Perelman}
G.~Perelman,
\newblock The entropy formula for the Ricci flow and its geometric
applications,
\newblock preprint, \href{https://arxiv.org/abs/math/0211159}
{arXiv:math/0211159} (2002).

\bibitem[VZ22]{VerdianiZiller}
L.~Verdiani and W.~Ziller,
\newblock Smoothness conditions in cohomogeneity one manifolds,
\newblock \emph{Transform. Groups} \textbf{27} (2022), 311--342.

\bibitem[Wan26a]{WangNegative}
Q.~S.~Wang,
\newblock Computer-assisted construction of
$\SU(2)$-invariant negative Einstein metrics,
\newblock \emph{Ann. Global Anal. Geom.} \textbf{69} (2026), no.~3,
Article 18.

\bibitem[Wan26b]{WangSpheres}
Q.~S.~Wang,
\newblock Doubly warped product Einstein metrics on spheres,
\newblock to appear in \emph{Comm. Anal. Geom.};
\href{https://arxiv.org/abs/2606.05519}{arXiv:2606.05519} (2026).

\bibitem[Win21]{WinkPotential}
M.~Wink,
\newblock Complete Ricci solitons via estimates on the soliton potential,
\newblock \emph{Int. Math. Res. Not.} \textbf{2021}, no.~6, 4487--4521.
\newblock \href{https://doi.org/10.1093/imrn/rnaa147}
{doi:10.1093/imrn/rnaa147}.

\bibitem[Win23]{WinkHopf}
M.~Wink,
\newblock Cohomogeneity one Ricci solitons from Hopf fibrations,
\newblock \emph{Comm. Anal. Geom.} \textbf{31} (2023), no.~3, 625--676.
\newblock \href{https://doi.org/10.4310/CAG.2023.v31.n3.a4}
{doi:10.4310/CAG.2023.v31.n3.a4}.

\bibitem[Zgl02]{ZgliczynskiLohner}
P.~Zgliczy\'nski,
\newblock $C^1$ Lohner algorithm,
\newblock \emph{Found. Comput. Math.} \textbf{2} (2002), 429--465,
\href{https://doi.org/10.1007/s102080010025}{doi:10.1007/s102080010025}.

\end{thebibliography}
\end{document}